          \documentclass{article}
\usepackage[paperwidth=7in, paperheight=10in, margin=.875in]{geometry}
\usepackage{subfig}
\usepackage{amsmath}
\usepackage{amssymb}
\usepackage{graphicx}
\usepackage{dcolumn}
\usepackage{mathtools}
\usepackage{amsmath,amsfonts}
\usepackage{bm}
\usepackage{amsmath}
\usepackage{amssymb}
\usepackage{algorithm2e}
\usepackage{color}
\usepackage{algorithmic}
\usepackage{float}

            \newtheorem{thm}{Theorem}[section]
          \newtheorem{prop}[thm]{Proposition}
          \newtheorem{lem}[thm]{Lemma}

          \newtheorem{ass}[thm]{Assumption}
\newcommand{\p}{\partial}
\newcommand{\f}{\frac}
\newcommand{\ds}{\displaystyle}
\newcommand{\E}{ {\mathbb{E}} }
\newcommand{\ola}{\overleftarrow}
\newcommand{\ora}{\overrightarrow}
\newcommand{\PP}{\mathbb{P}}
\newcommand{\F}{\mathcal{F}}
\newcommand{\R}{\mathbb{R}}
\newcommand{\be}{\begin{equation}}
\newcommand{\ee}{\end{equation}}
\newcommand{\M}{\mathcal{M}}
\newcommand{\sS}{\mathcal{S}}
\newcommand{\C}{\mathcal{C}}
\renewcommand{\L}{\mathcal{L}}
\renewcommand{\d}{\mathrm{d}}
\def\la{\left\langle}
\def\ra{\right\rangle}
 \date{}
           
\begin{document}
\title{Forward Backward Doubly Stochastic Differential Equations and the Optimal Filtering of Diffusion Processes \thanks{}}
\author{Feng Bao \thanks{ Department of Mathematics, The University of Tennessee at Chattanooga, Chattanooga, Tennessee, 37403 \ ({\tt feng-bao@utc.edu}).}
        \and Yanzhao Cao \thanks{Department of Mathematics and Statistics, Auburn University, Auburn, Alabama,
        36849 \ ({\tt yzc0009@auburn.edu}).}
        \and Xiaoying Han \thanks{Department of Mathematics and Statistics, Auburn University, Auburn, Alabama,
        36849 \ ({\tt xzh0003@auburn.edu}).}
        }

\pagestyle{myheadings} 

\markboth{Optimal Filtering of Diffusion Processes}{Bao, Cao \& Han} \maketitle

\begin{abstract}
The connection between forward backward doubly stochastic differential equations and the optimal filtering problem is established without using the Zakai's equation.    The solutions of forward backward doubly stochastic differential equations are expressed in terms of conditional law of a partially observed Markov diffusion process.   It then follows that the adjoint time-inverse forward backward doubly stochastic differential equations governs the evolution of the unnormalized filtering density in the optimal filtering problem.
\end{abstract}


\section{Introduction}

The goal of this work is to study the state of a noise-perturbed dynamical system,  $U_t$, given noisy observation on the dynamics, $V_t$.   This suggests the optimal filtering problem of determining the conditional probability of $U_t$, given an observed path $\{ V_s: 0 \leq s \leq t \}$.   The pioneer work of optimal filtering problems was considered by  Kallianpur and Striebel \cite{Kallianpur} and Zakai \cite{zakai}.   In particular, the Kallianpur-Striebel formula provides a continuous time framework of the optimal filtering that considers the conditional probability density function (PDF) of the state as the solution of a nonlinear stochastic partial differential equation (SPDE);  and the approach proposed by Zakai leads to a linear stochastic integro-differential parabolic equation, referred to as the Zakai's  equation.   Under strong regularity conditions it can be shown that the solution of the Zakai's equation represents an unnormalized conditional density of the state process.
Fundamental research of the optimal filtering problem was also conducted  by Kalman and Bucy \cite{Bucy, Kalman}, Kushner and Pardoux \cite{Kushner, Pardoux_SPDE}, Shiryaev \cite{Shiryaev} and Stratonovich  \cite{Stratonovich}, among other extensive studies on discrete nonlinear filter solver (see \cite{Doucet_PF, Dun-EKF, Gobet-Zakai, particle-filter, Julier-EKF}). 
\smallskip

The advantage of solving the optimal filter problems with SPDEs such as the Zakai equation is that it provides the ``exact" solution for the conditional density of $U_t$ given   $\{ V_s\}_{0 \leq s \leq t }$. However, it has not been considered as an efficient method by the science and engineering community because of its slow convergence and high complexity.   Instead of dealing with SPDEs,   the unnormalized density function can also be studied through a system of  stochastic (ordinary) differential equations (SDEs). Such a system consists of two SDEs,  one standard SDE and one backward doubly stochastic differential equation (BDSDE), and is referred to as a system of forward backward doubly stochastic differential equations (FBDSDEs).  
The FBDSDE system was first studied by  Pardoux and Peng in \cite{PP1994},  where  the equivalence between FBDSDEs and certain parabolic type SPDEs was established.   Our recent work  \cite{BC2014, Bao_Full,  Bao, Bao_Semi} indicates  that solving optimal filtering problems with FBDSDE  systems  can be far less costly than that with SPDEs and more accurate than both SPDEs and discrete filter methods such as particle filter methods. 

\smallskip
In this paper, we establish a direct link between the optimal filtering problem and a FBDSDE system. First we provide a FBDSDE version of Feynman-Kac formula for the optimal filter problem and obtain the  adjoint  of this system.  To the best of our knowledge, similar  results have been obtained before.  As a consequence, we show this adjoint, which is a  a time-inverse FBDSDE system,  provides a solution for the unnormalized condition density of the optimal filter problem. 

\smallskip
The rest of this paper is organized as follows.  In Section \ref{sec:pre} we present the  mathematical formulation of the optimal filtering problem and provide a brief introduction of  FBDSDEs.  In Section \ref{sec:main} we establish the connection between the FBDSDEs and the unnormalized conditional density function.  Some closing remarks will be given in Section \ref{sec:close}.

\section{Preliminaries}\label{sec:pre}

In this section, we present the  mathematical formulation of the optimal filtering problem and provide a brief introduction of  FBDSDEs.   

\smallskip
Let $(\Omega, \F, \PP)$ be a probability space, and let $T > 0$ be fixed throughout the paper.    Let $\{W_t\}_{0 \leq t \leq T}$ and $\{B_t\}_{0 \leq t \leq T}$ be two mutually independent standard Brownian motions defined on  $(\Omega, \F, \PP)$, with values in $\mathbb{R}^d$ and $\mathbb{R}^l$, respectively.  
Denote by $\mathcal{N}$ the class of $\PP$-null sets of $\F$.  For each $t \in [0, T]$ and any process $\eta_t$,  let 
$$\mathcal{F}_{s,t}^{\eta}: = \sigma \{\eta_r - \eta_s: s\leq r \leq t\} \vee  \mathcal{N}$$   be the $\sigma$-field generated by $\{\eta_r - \eta_s\}_{s\leq r \leq t}$ and write  $\mathcal{F}_t^{\eta} = \mathcal{F}_{0,t}^{\eta}$.


\subsection{The optimal filtering problem}

Consider the following stochastic differential system on the probability space $(\Omega, \mathcal{F}, \PP)$
\begin{equation}\label{NLF}
\left\{
\begin{aligned}
\d U_t =& b_t(U_t) \d t + \rho_{t} \d W_t + \tilde{\rho}_{t} \d B_t, \\[1.2ex]
\d V_t =& h(U_t) \d t + \d B_t,
\end{aligned} \right.
\end{equation}
where $\{U_t \in \mathbb{R}^d: t \geq 0\}$  is the ``state process" that describes the state of a dynamical system and $\{V_t \in \mathbb{R}^l: t \geq 0\}$ is the ``measurement process'' which is the noise perturbed observations of the state $U_t$.     Given an initial state $U_0$ with probability distribution $p_0(u)$  independent of $W_t$ and $B_t$, the goal of the optimal filtering problem is to obtain the best estimate of $\phi(U_t)$ as the conditional expectation with respect to the measurement $\{V_r\}_{ 0 \leq s \leq t}$, where $\phi$ is a given test function.

\smallskip

Denote by $\F^V_t: = \sigma \{V_r : 0 \leq r \leq t \}$ the $\sigma$-field generated by the measurement process from time $0$ to $t$ and denote by $\mathcal{M}_t$   the space of all $\mathcal{F}_t^V$-measurable and square integrable random variables at time $t$.  The optimal filtering problem can be formulated mathematically as to find the conditional expectation
 \begin{equation*}
\E \left[\phi(U_t) \big| \mathcal{F}_t^V \right] = \inf\left\{ \E\left[ |\phi(U_t) - \psi_t\right|^2]:  \psi_t \in \mathcal{M}_t \right\}.
\end{equation*}
According to \cite{Kallianpur1968, Kallianpur1969},  the {\it optimal} filter is given by
\begin{equation}\label{eq:Kallianpur} 
\E \left[\phi(U_t) \big| \mathcal{F}_t^V \right] = \f{\displaystyle\int_{\mathbb{R}^d} \phi(u) p_t du}{ \displaystyle\int_{\mathbb{R}^d} p_t du },
\end{equation}
where $p_t$ is the {\it unnormalized filtering density}.  \eqref{eq:Kallianpur} is the well known Kallianpur--Striebel formula.

\smallskip
Define  
$$Q_t^s := \exp \left\{ \int_s^t h(U_r) d U_r - \f{1}{2} \int_s^t | h(U_r) |^2 dr \right\}. $$
When $s=0$ we denote  $Q_t^{0}$ as $Q_t$ in short.  Let $\tilde{\PP}$ be the probability measure induced on the space $(\Omega, \mathcal{F})$ such that
\begin{equation}\label{eq:ptilde}
\f{d \PP}{d\tilde{\PP}}\bigg|_{\mathcal{F}_t^V} = Q_t.
\end{equation}
Then according to the Cameron-Martin theorem  the probability measures $\PP$ and $\tilde{\PP}$ are equivalent  when the Novikov condition is satisfied \cite{Girsanov}.  
Moreover, it is straightforward to verify that (see \cite{Oksendal2003}, Lemma 8.6.2)
\begin{equation}\label{eq:pro-transform}
\E\left[\phi(U_t) \big| \mathcal{F}_t^V\right]=\frac{ \tilde{\E}\left[ \phi(U_t) Q_t \big| \mathcal{F}_t^V \right]}{ \tilde{\E}\left[ Q_t \big| \mathcal{F}_t^V \right]}.
\end{equation}
where  $\tilde{\E}$  denotes  the expectation with respect to  $\tilde{\PP}$.

\subsection{Forward backward doubly stochastic differential equations}

For each $t \in [0, T]$, define
$$
\mathcal{F}_t := \mathcal{F}_t^{W} \vee \mathcal{F}_{t,T}^B.
$$
Then the collection $\{\mathcal{F}_t: t\in [0,T]\}$ is neither increasing nor decreasing, and thus   does not constitute a filtration \cite{PP1994}.     For any positive integer $n \in \mathbb{N}$,   denote by $\M^2(0,T; \mathbb{R}^n)$ the set of $\mathbb{R}^n$-valued jointly measurable random processes $\{\psi_t: t \in [0, T] \}$ such that $\psi_t$ is $\mathcal{F}_t$ measurable for a.e. $t \in [0, T]$ and satisfies
$$\E \int_0^T |\psi_t|^2 dt < \infty. $$
Similarly,   denote by $\sS^2([0,T];\mathbb{R}^n)$ the set of   continuous $\mathbb{R}^n$-valued random processes $\{\psi_t: t \in [0, T] \}$ such that $\psi_t$ is $\mathcal{F}_t$ measurable for any $t \in [0, T]$ and satisfies
$$ \E  \sup_{0 \leq t \leq T}|\psi_t|^2  	 < \infty. $$

\smallskip
We next provide a brief introduction of forward backward doubly stochastic differential equations (FBDSDEs), summarized from \cite{PP1994}.   

Given $\tau \geq 0$, $x \in \R^d$ and $\varphi \in \L^2(\Omega, \F_T, \PP)$,  a system of forward backward doubly stochastic differential equations (FBDSDEs) can be formulated as
\begin{eqnarray*}
\ds \d X_t &=&   b(X_t) \d t + \sigma(X_t) \d W_t, \quad \tau  \leq t \leq T,  \\[1.2ex]
 \ds -\d Y_t &=& \ds    f(t, X_t, Y_t, Z_t) \d t   +  \ds g(t, X_t, Y_t,  Z_t) \d \ola{B}_{t} - Z_t \d W_{t},\quad \tau  \leq t \leq T,  \\[1.2ex]
 \ds X_{\tau} &=&  x, \qquad Y_T   =  \varphi(X_T),  
\end{eqnarray*}
or, in the integral equation form,  for any $t \in [\tau, T]$, 
 \begin{eqnarray}
    X_t &=&  x + \int_\tau^t b(X_s) \d s + \int_\tau^t \sigma(X_s) \d W_s, \label{FBDSDEs:eq1}\\[1.2ex]
    Y_t &=&  \varphi(X_T) + \int^T_t   f(s, X_s, Y_s, Z_s)\d s  +  \int^T_t g(s, X_s, Y_s,  Z_s) \d \ola{B}_{s}  -  \int^T_t Z_s \d W_{s}.
\label{FBDSDEs:eq2}
\end{eqnarray}
    Notice that equation \eqref{FBDSDEs:eq1} is a standard forward SDE with a standard forward  It\^o integral and equation \eqref{FBDSDEs:eq2} is a backward doubly stochastic differential equation (BDSDE) involving the backward It\^o integral $\int \cdot \d \ola{B}_{s}$ (see \cite{Two_sided} for details on the  two types of integrals).

\smallskip
Let the mappings $f: [0,T] \times \mathbb{R}^d \times \mathbb{R}^k \times \mathbb{R}^{k \times d} \rightarrow \mathbb{R}^k $
and
$g:   [0,T]  \times \mathbb{R}^d \times \mathbb{R}^k \times \mathbb{R}^{k \times d} \rightarrow \mathbb{R}^{k\times l}$
be jointly measurable and for any $(y,z) \in   \mathbb{R}^k \times \mathbb{R}^{k\times d}$,
$$f(\cdot, \cdot, y, z)\in \M^2(0,T; \mathbb{R}^k),
\qquad g(\cdot, \cdot, y, z)\in \M^2(0,T; \mathbb{R}^{k\times l}). $$     Denote by $| \cdot |$ the Euclidean norm of a vector and by $\|A \|: = \sqrt{\mathrm{Tr} ( A A^*)}$ the norm of a matrix $A$.   The existence and uniqueness of solutions, moment estimates for the solutions, and the regularity of solutions to Equation \eqref{FBDSDEs:eq2} rely on one or more of the following assumptions.

\smallskip
\begin{ass}\label{Lipschitz}
$f$ and $g$ satisfy the Lipschitz condition:  there exist   constants $c>0$ and $0 < \bar{c} < 1$ such that for any   $(t, x) \in  [0, T] \times  \R^d$,  $y_1, y_2 \in \mathbb{R}^k $ and $z_1, z_2 \in \mathbb{R}^{k \times d}$,
\begin{eqnarray*}
|f(t, x, y_1, z_1) - f(t, x, y_2, z_2)|^2 &\leq& c  ( |y_1 - y_2|^2 + \|z_1 - z_2\|^2) ,
\\[1.2ex]
 \|g(t, x, y_1, z_1 ) - g(t, x, y_2, z_2 )\|^2 &\leq& c |y_1 - y_2|^2 + \bar{c} \|z_1 - z_2\|^2.
\end{eqnarray*}
\end{ass}
\begin{ass}\label{g:bounded}
There exists $c >0$ such that for all $(t, x, y, z) \in [0, T] \times \R^d \times \mathbb{R}^k \times \mathbb{R}^{k \times d}$, 
$$g g^{\ast}(t, x, y, z) \leq z z^{\ast} + c (\|g(t, x, 0, 0)\|^2 + |y|^2) I.$$
\end{ass}
\begin{ass}\label{gprime:bound}
For any $(t, x, y, z) \in [0, T] \times \R^d \times \mathbb{R}^k \times \mathbb{R}^{k \times d}$ and $ \theta \in \R^{k \times d}$
$$\f{\p g}{\p z} (t, x, y, z) \theta \theta^* \left(\f{\p g}{\p z} (t, x, y, z)\right)^* \leq \theta \theta^*.$$
\end{ass}

The following results are due to Pardoux and Peng \cite{PP1994}.

\begin{prop}\label{exist-uniq} Under Assumption \ref{Lipschitz}, the BDSDE \eqref{FBDSDEs:eq2} admits a unique solution 
$$(Y, Z)    \in  \sS^2([0,T]; \mathbb{R}^k) \times \M^2(0,T; \mathbb{R}^{k\times d}). $$
\end{prop}

\begin{prop}\label{Y:bounded} 
Let Assumptions \ref{Lipschitz} and \ref{g:bounded} hold, then the solution  of the BDSDE \eqref{FBDSDEs:eq2} satisfies  $$ \E \ds \sup_{0 \leq t \leq T } |Y_t |^2  < \infty.  $$
\end{prop}

For any positive integer $k$, denote by $\C^k_{l,b}$ the collection of $\C^k$ functions with bounded partial derivatives of all orders less than or equal to $k$, and denote by $\C^k_p$ the collection of $\C^k$ functions with partial derivatives of all orders less than or equal to $k$ which grow at most like a polynomial function of $x$ as $x\rightarrow \infty$.    It is well known that given $b \in \C^3_{l, b} (\R^d, \R^d)$ and $\sigma \in \C^3_{l,b}(\R^d, \R^{d \times d})$,  for each $(\tau, x) \in [0, T] \times \R^d$, the SDE \eqref{FBDSDEs:eq1} has a unique strong solution, denoted as $X_t^{\tau, x}$.   Consequently denote by $(Y_t^{\tau, x}, Z_t^{\tau, x})$ the unique solution to the BDSDE
    \be  \label{FBDSDEs:eq3} Y_t =  \varphi(X_T^{\tau, x}) + \int^T_t   f(s, X^{\tau, x}_s, Y_s, Z_s)\d s  +  \int^T_t g(s, X^{\tau, x}_s, Y_s,  Z_s) \d \ola{B}_{s}  -  \int^T_t Z_s \d W_{s}.  \ee
    
\begin{prop}\label{reg}
Let $\varphi \in \C^3_p(\R^d; \R^k)$.   Under Assumptions \ref{Lipschitz} -- \ref{gprime:bound},    the random field $\left\{Y^{\tau, x}_\tau: \tau \in [0, T], x \in \R^d \right\}$ admits a continuous version such that for any $\tau \in [0, T]$, $x \mapsto Y^{\tau, x}_\tau$ is of class $\C^2$ a.s..
\end{prop}

 \smallskip
 
The following regularity result   can be obtained by using  standard techniques of SDEs, FBSDEs and BDSDEs (see Proposition 1 in \cite{Bao}) \begin{lem}\label{Y:regularity}
In addition to the Assumption \ref{Lipschitz}, assume that $f, g \in \C_{l, b}^{1}$.  Then  the solution $(Y_t^{\tau, x}, Z_t^{\tau, x})$ to the BDSDE \eqref{FBDSDEs:eq3} satisfies
$$ \E\left[ \left(Y_t^{\tau, x} - Y_{\tau}^{\tau, x}\right)^2 \right] \leq C(t - \tau), \quad \E\left[ \left(Z_t^{\tau, x} - Z_{\tau}^{\tau, x}\right)^2 \right] \leq C(t - \tau), \quad 0 \leq \tau \leq t \leq T, $$
where $C$ is a positive constant independent of $\tau$ and $t$.
\end{lem}

\smallskip
Note that with the convention above, the unique solution to the FBDSDE system \eqref{FBDSDEs:eq1} -- \eqref{FBDSDEs:eq2} can be written as $(X^{\tau, x}_t, Y^{\tau, x}_t, Z^{\tau, x}_t)$.   
Denote
$$\nabla X_t^{\tau, x} := \f{\p X_t^{\tau, x}}{\p x}, \qquad  \nabla Y_t^{\tau, x}:= \f{\p Y_t^{\tau, x}}{\p x}, \qquad  \nabla Z_t^{\tau, x}:= \f{\p Z_t^{\tau, x}}{\p x}.$$   Then $(\nabla Y_t^{\tau, x}, \nabla Z_t^{\tau, x})$ is the unique solution to  variational form of the BDSDE \eqref{FBDSDEs:eq2} (see \cite{PP1994})
\begin{equation*}\label{Y:variance}
\begin{aligned}
 \nabla Y_t^{\tau, x}  &=  \varphi^\prime(X^{\tau, x}_T) \nabla X_T^{\tau, x} + \int^T_t \left( \f{\p f}{\p x} \nabla X_s^{\tau, x} + \f{\p f}{\p Y}\nabla Y_s^{\tau, x} + \f{\p f}{\p Z}\nabla Z_s^{\tau, x} \right) \d s  \\
& \quad+ \int^T_t \left( \f{\p g}{\p x} \nabla X^{\tau, x}_s + \f{\p g}{\p Y} \nabla Y_s^{\tau, x} + \f{\p g}{\p Z}\nabla Z_s^{\tau, x} \right) \d \ola{B}_s 
 - \int^T_t  \nabla Z_s^{\tau, x} \d W_s.
\end{aligned} 
\end{equation*}
In addition, the random field $\left\{Z^{\tau, x}_t :  t \in [\tau, T], x \in \R^d\right\}$ has an a.s. continuous version 
\begin{equation} \label{Z=dY}
Z_t^{\tau,x} = \nabla Y_t^{\tau, x} (\nabla X_t^{\tau, x})^{-1} \sigma(X_t^{\tau, x}),  \qquad Z_{\tau}^{\tau,x} = \nabla Y_{\tau}^{\tau, x} \sigma(x) .
\end{equation}

\smallskip

The following Lemma follows directly from Lemma \ref{Y:regularity} and  Proposition \ref{Y:bounded}.
\begin{lem}\label{dY:regularity}
Assume that $b \in \C_{l, b}^{2}$, $f \in \C_{l, b}^{2}$, $g \in \C_{l, b}^{2}$ and $\varphi \in \C_{l, b}^{2}$.  Then there exists $C > 0$ such that
$$ \E[ (\nabla Y_t^{\tau, x} -  \nabla Y_t^{\tau, x})^2 ] \leq C(t - \tau), \quad \E[ (\nabla Z_t^{\tau, x} - \nabla Z_{\tau}^{\tau, x})^2 ] \leq C(t - \tau), \quad 0 \leq \tau \leq t \leq T.$$
Moreover, 
  $$\E \ds \sup_{0 \leq t \leq T } |\nabla Y_t^{\tau, x} |^2  < \infty.$$
\end{lem}

\section{FBDSDEs and Optimal Filtering}\label{sec:main}

In this section, we establish the connection between  the optimal filtering problem and a FBDSDE system.  
In particular, we will first prove a Feynman-Kac formula in the filtering context. Then we present the adjoint relationship between standard FBDSDEs and time-inverse FBDSDEs.  In the end we will show that the solution of a time-inverse FBDSDE  is  the unnormalized filtering density sought in the optimal filtering problem. For simplicity of exposition, we only discuss the one dimensional case with $d=1$ and $l=1$.   The same method can also be applied to multi-dimensional cases  with more complicated calculations.

\subsection{Feynman-Kac type formula for optimal filtering}

For $\tau \in [0, T]$ and $x \in \R^d$, consider the following FBDSDE system on the probability space $(\Omega, \F, \tilde{\PP})$  
\begin{equation}\label{Backward:FBDSDEs}
\left\{
\begin{aligned}
\d X_t &=  b_t(X_t) \d t + \sigma_t \d W_t,   \quad \tau \leq t \leq T & \text{ (SDE) } \\[1.2ex]
- \d Y_t &=- Z_t  \d W_t + \left( h(X_t )  Y_t  + \f{\tilde{\rho}_t}{\sigma_t} Z_t  \right) \d \ola{V}_t, \quad \tau \leq t \leq T & \text{(BDSDE)}\\[1.2ex]
X_{\tau} &= x,  \qquad  Y_T   =  \phi(X_T), 
\end{aligned}\right.
\end{equation}
where $\sigma_t^2 = \rho_t^2 + \tilde{\rho}_t^2$, and $b$, $\rho$, $\tilde{\rho}$ , $h$ are the functions  appeared in the optimal filtering problem \eqref{NLF}.     Here $W_t$ is the same Brownian motion as in the nonlinear filtering problem \eqref{NLF},  while  $V_t$ is the measurement process  which becomes a standard Brownian motion independent of $W_t$ under the induced probability measure $\tilde{\PP}$ defined by \eqref{eq:ptilde}.    Then $X_t$ is a $\mathcal{F}^{W}_{t}$ adaptive stochastic process and the pair $( Y_t,  Z_t )$ is adaptive to  $\mathcal{F}^{W}_{t} \vee \mathcal{F}^{V}_{t, T}$.    For any  single-variable function $F=F(x)$, denote $F^{\prime} : = \f{d F}{d x}$  and $F^{\prime\prime} = \f{d^2 F}{d x^2}$.

\smallskip
\begin{lem}\label{h:regularity}
Assume that $ b_t$ and $\sigma_t$ are bounded and  $h \in \C_b^{2}(\mathbb{R}; \R)$.   Then for any $0 \leq s \leq t \leq T$, there exists a positive constant $C$   independent of $s$ and $t$ such that
\begin{equation}\label{eq:h-regularity}
\tilde{\E}[ ( h(X_{t}) - h(X_{s}) )^2 | \mathcal{F}_{t,T}^{V}] \leq C (t - s).
\end{equation}
\end{lem}

\noindent \textit{Proof.}   The application of   It\^o's formula to $ h(X_{t})$ results in
$$
h(X_{t}) = h(X_s) + \int_{s}^{t} \left( b_r (X_{r}) h^{\prime}(X_{r})  + \f{\sigma_r^2}{2} h^{\prime \prime}(X_{r}) \right) \d r + \int_{s}^{t} \sigma_r  h^{ \prime}(X_{r})  \d W_r,
$$
and hence
\be \label{eq:hsquare}
\left( h(X_{t}) - h(X_{s}) \right)^2   = \left( \int_{s}^{t} \left( b_r (X_{r}) h^{\prime}(X_{r})  + \f{\sigma_r^2}{2} h^{\prime \prime}(X_{r}) \right) \d r + \int_{s}^{t} \sigma_r  h^{ \prime}(X_{r})  \d W_r \right)^2. 
\ee
Taking expectation $\tilde{\E}$ of the above gives
$$
\tilde{\E}\big[ ( h(X_{t}) - h(X_{s}) )^2 \big] = \tilde{\E}\left[\Big( \int_{s}^{t} \big( b_r (X_{r}) h^{\prime}(X_{r})  + \f{\sigma_r^2}{2} h^{\prime \prime}(X_{r}) \big) \d r \Big)^2 \right] + \tilde{\E}\left[ \int_{s}^{t}  ( \sigma_r  h^{ \prime}(X_{r}) \big )^2  \d r \right].$$
The inequality \eqref{eq:h-regularity} then follows immediately from the assumptions of the lemma.  \hfill $\Box$

\smallskip
With Proposition \ref{Y:bounded}, and Lemmas \ref{Y:regularity}  and \ref{h:regularity},  we establish the following    Feynman-Kac formula in the optimal filtering context.

\begin{thm}\label{thm:Feynman-Kac}
Assume that $\phi$ is bounded,  $b_t, \rho_t, \tilde{\rho}_t \in \C_{l, b} $ and $h \in \C_{l, b}^{2}(\mathbb{R})$.  Then, $\forall \tau \in [0, T]$ and $x \in \mathbb{R}^d$ the following equality holds a.s.
\begin{equation}\label{Feynman-Kac-Filter}
Y_\tau^{\tau, x} = \E_{\tau}^{x}[\phi(U_T) Q_T^\tau ],
\end{equation} 
where $ \E_{\tau}^{x}[\cdot] : = \tilde{\E}[\cdot | \mathcal{F}_{\tau,T}^{V}, U_\tau = x ]$. 
\end{thm}

\smallskip
\noindent \textit{Proof.}   We prove  the statement \eqref{Feynman-Kac-Filter} for $\tau = 0$ only, the general case follows from the 
$\tau = 0$  case trivially.  First it is straightforward to  verify that under assumptions in Theorem \ref{thm:Feynman-Kac}, all the assumptions of Proposition \ref{Y:bounded}, and Lemmas \ref{Y:regularity}  and \ref{h:regularity} are fulfilled.   Since  $Y_\tau^{\tau,x}$ and $Z_\tau^{\tau,x}$ are functions of $x$, we write
 $Y_\tau^{\tau,x} = Y_\tau(x)$ and $Z_\tau^{\tau,x} = Z_\tau(x) $ in the sequel.  \smallskip

 Let $0 = t_0 < t_1 < t_2 \cdots < t_N = T$ be an equidistant temporal partition with $t_{n+1} - t_n = {T}/{N} := \Delta t$ and define $$\Delta_n = \E_{0}^{x}[Q_{t_{n+1}} Y_{t_{n+1}}(U_{t_{n+1}}) - Q_{t_n} Y_{t_n}(U_{t_n})].$$  It follows immediately that
$$  \E_{0}^{x}[\phi(U_T) Q_T  - Y_0(x)]  = \sum_{n=0}^{N-1} \Delta_n. $$
 Denote $\tilde{\PP}_x:= \tilde{\PP}(\cdot | U_0 = x)$.  To prove  \eqref{Feynman-Kac-Filter} it suffices to verify that 
$$\sum_{n = 0}^{N-1} \Delta_n \stackrel{N \to  \infty }{\longrightarrow} 0 \ \text{in} \ \L^1(\Omega, \tilde{\PP}_x).$$

For each $n \geq 0$, let $U_{t_n}$ be the solution of the state for \eqref{NLF} at time step $t_n$ and consider the FBDSDEs system \eqref{Backward:FBDSDEs} on $[t_n, t_{n+1}]$ with initial condition $U_{t_n}$:
\begin{equation}\label{FBDSDEs:interval}
\left\{
\begin{aligned}
\d \hat{X}_t & = b_t (\hat{X}_t) \d t + \sigma_t \d W_t, \\
- \d Y_t &= - Z_t \d W_t + \left( h(\hat{X}_t) Y_t  + \f{\tilde{\rho}_t}{\sigma_t} Z_t \right) \d \ola{V}_t, \\
\hat{X}_{t_n} &= U_{t_n}, \hspace{1em}
Y_{t_{n+1}} = Y_{t_{n+1}}(\hat{X}_{t_{n+1}}) .
\end{aligned} \right.
\end{equation}
From the definition of the state process $U_t$ in \eqref{NLF} and the SDE $\hat{X}_t$ in \eqref{FBDSDEs:interval},  we have the relation between $U_{t_{n+1}}$ and $\hat{X}_{t_{n+1}}$: \begin{equation*}
U_{t_{n+1}} = \hat{X}_{t_{n+1}} + \int_{t_n}^{t_{n+1}} \rho_s \d W_s  -  \int_{t_n}^{t_{n+1}} \sigma_s(\hat{X}_s) \d W_s + \int_{t_n}^{t_{n+1}} \tilde{\rho}_s \left( \d V_s -  h(U_s) \d s \right) + R_X^{n+1}
\end{equation*}
where $$R_X^{n+1} =  \int_{t_n}^{t_{n+1}} b_s (U_s) ds - \int_{t_n}^{t_{n+1}} b_s (\hat{X}_s) ds .$$

To simplify presentation, for any process $\psi_t$ we write  $\hat{\psi}_t := \psi_t (\hat{X}_t)$  throughout the rest of this proof.    Let $\eta_{n+1}  =  U_{t_{n+1}} - \hat{X}_{t_{n+1}}$ Then from the above we have that 
\begin{equation}\label{def:eta}
\eta_{n+1}  = \int_{t_n}^{t_{n+1}} \rho_s \d W_s  -  \int_{t_n}^{t_{n+1}} \sigma_s \d W_s + \int_{t_n}^{t_{n+1}} \tilde{\rho}_s \left( \d V_s -  h(U_s) \d s \right) + R_X^{n+1}
\end{equation}
Applying the Taylor expansion to $Y_{t_{n+1}}$ we have that 
\begin{equation}\label{Taylor:Y:X}
Y_{t_{n+1}}(U_{t_{n+1}})   = \hat{Y}_{t_{n+1}} + \hat{Y}^{'}_{t_{n+1}} \cdot \eta_{n+1} + \f{1}{2} \hat{Y}^{''}_{t_{n+1}} \cdot ( \eta_{n+1} )^2 +  \xi_{n+1},
\end{equation}
where  $\xi_{n+1}$ is the Taylor remainder such that  $\E_0^x[(\xi_{n+1})^2] \leq C (\Delta t)^{3}$.    Then for each $n = 0, 1, 2, \cdots, N-1$,
\begin{equation}\label{Delta_n:step1}
\begin{aligned}
\Delta_n &=  \E_{0}^{x} \left[ Q_{t_{n+1}}Y_{t_{n+1}}(U_{t_{n+1}}) - Q_{t_{n}} \hat{Y}_{t_{n+1}} + Q_{t_{n}} \hat{Y}_{t_{n+1}} - Q_{t_{n}} Y_{t_n}(U_{t_n})\right] \\[1.2ex]
 & =  \underbrace{\E_{0}^{x} \left[ \left( Q_{t_{n+1}}  - Q_{t_{n}} \right) \hat{Y}_{t_{n+1}}\right]}_{({\bf i})} + \underbrace{\E_{0}^{x}\left[ Q_{t_{n}} \left( \hat{Y}_{t_{n+1}} - Y_{t_n}(U_{t_n}) \right)\right]}_{({\bf ii})} \\
& \qquad + \underbrace{\E_{0}^{x}\left[ Q_{t_{n+1}} \left( \hat{Y}^{'}_{t_{n+1}} \eta_{n+1} + \f{1}{2} \hat{Y}^{''}_{t_{n+1}}\cdot ( \eta_{n+1} )^2 +  \xi_{n+1} \right) \right]}_{({\bf iii})}.
\end{aligned}
\end{equation}
We next estimate terms ({\bf i}), ({\bf ii}) and ({\bf iii}) in \eqref{Delta_n:step1} one by one.

\smallskip
\noindent ({\bf i})   Write $h_t = h(U_t)$ and $\hat{h}_t = h(\hat{X}_t)$, and apply Ito's formula to $Q_{t_{n}}$ we obtain
\begin{equation}\label{dQ:term1}
\begin{aligned}
 \E_{0}^{x} \left[ ( Q_{t_{n+1}} - Q_{t_{n}} ) \hat{Y}_{t_{n+1}}\right]  = & \E_{0}^{x} \left[\int_{t_n}^{t_{n+1}} h_s Q_s \d V_s \hat{Y}_{t_{n+1}} \right] \\
= &  \E_{0}^{x} \left[\int_{t_n}^{t_{n+1}} \hat{h}_s Q_s \d V_s \hat{Y}_{t_{n+1}}\right] + \E_{0}^{x} \left[\int_{t_n}^{t_{n+1}} (h_s - \hat{h}_s) Q_s \d V_s \hat{Y}_{t_{n+1}} \right].
\end{aligned}
\end{equation}
Applying It\^o formula to function $h$ yields
$$h_s - \hat{h}_s = h^{'}(U_{t_n})\left( \int_{t_n}^{s} \rho_r \d W_r + \int_{t_n}^{s} \tilde{\rho}_r \d V_r - \int_{t_n}^{s} \sigma_r \d W_r \right) + \mathcal{O}(\Delta t), $$
and consequently with  $h^{'}_{t_n} := h^{'}(U_{t_n})$ we have
\begin{equation}\label{h-h_bar}
\begin{aligned}
& \E_{0}^{x} \left[\int_{t_n}^{t_{n+1}} (h_s - \hat{h}_s) Q_s \d V_s \hat{Y}_{t_{n+1}}\right] \\
= & \E_{0}^{x} \left[  h^{'}_{t_n} Q_{t_n} Y_{t_{n}}( U_{t_{n}})\int_{t_n}^{t_{n+1}} \d V_s  \left( \int_{t_n}^{s} \rho_r \d W_r  - \int_{t_n}^{s} \sigma_r \d W_r \right)  \right] \\
& \quad  + \E_{0}^{x} \left[  h^{'}_{t_n} \left(  Q_s - Q_{t_n} \right) \left( \hat{Y}_{t_{n+1}} - Y_{t_{n}}( U_{t_{n}})  \right) \int_{t_n}^{t_{n+1}} \d V_s  \left( \int_{t_n}^{s} \rho_r \d W_r  - \int_{t_n}^{s} \sigma_r \d W_r \right)  \right] \\
& \quad +  \E_{0}^{x} \left[h^{'}_{t_n} Q_{t_n} \hat{Y}_{t_{n+1}}  \int_{t_n}^{t_{n+1}} \int_{t_n}^{s} \tilde{\rho}_r \d V_r  \d V_s \right]  + \mathcal{O}\left((\Delta t)^{\f{3}{2}}\right).
\end{aligned}
\end{equation}

\smallskip
First noting that $h^{'}_{t_n} Q_{t_n} Y_{t_{n}}( U_{t_{n}})\int_{t_n}^{t_{n+1}} \d V_s $ is independent of $\int_{t_n}^{s} \rho_r \d W_r  - \int_{t_n}^{s} \sigma_r \d W_r$, we have
\be \label{i:est1}\E_{0}^{x} \left[  h^{'}_{t_n} Q_{t_n} Y_{t_{n}}( U_{t_{n}})\int_{t_n}^{t_{n+1}} \d V_s  \left( \int_{t_n}^{s} \rho_r \d W_r  - \int_{t_n}^{s} \sigma_r \d W_r \right)  \right] = 0. \ee
Second, it's straightforward to verify that 
\be \label{i:est2}  \E_{0}^{x} \left[  h^{'}_{t_n} \left(  Q_s - Q_{t_n} \right) \left( \hat{Y}_{t_{n+1}} - Y_{t_{n}}( U_{t_{n}})  \right) \int_{t_n}^{t_{n+1}} \d V_s  \left( \int_{t_n}^{s} \rho_r \d W_r  - \int_{t_n}^{s} \sigma_r \d W_r \right)  \right] \sim \mathcal{O}\left((\Delta t)^{\f{3}{2}}\right).\ee
Putting \eqref{i:est1} and \eqref{i:est2} in \eqref{h-h_bar}, it follows from the regularity condition of $\tilde{\rho}_r$ that 
$$
\begin{aligned} 
& \E_{0}^{x} \left[\int_{t_n}^{t_{n+1}} (h_s - \hat{h}_s) Q_s \d V_s \hat{Y}_{t_{n+1}}\right]  = \E_{0}^{x} \left[h^{'}_{t_n} Q_{t_n} \hat{Y}_{t_{n+1}} \tilde{\rho}_{t_n} \int_{t_n}^{t_{n+1}} \int_{t_n}^{s}  \d V_r  \d V_s \right]  + \mathcal{O}\left((\Delta t)^{\f{3}{2}}\right) .
\end{aligned}
$$

\smallskip
Define
\be \label{def:nu} \nu_n: = h^{'}_{t_n} Q_{t_n} \hat{Y}_{t_{n+1}} \tilde{\rho}_{t_n} \int_{t_n}^{t_{n+1}} \int_{t_n}^{s}  \d V_r  \d V_s. \ee   Then by using
 the facts $\int_{t_n}^{t_{n+1}} \int_{t_n}^{s}  \d V_r  \d V_s = \f{1}{2} \big( (V_{t_{n+1}} - V_{t_n})^2 - \Delta t \big)$ and $h^{'}_{t_n} Q_{t_n} \hat{Y}_{t_{n+1}} \tilde{\rho}_{t_n}$ is independent of $\f{1}{2} \big( (V_{t_{n+1}} - V_{t_n})^2 - \Delta t \big)$ we have
\begin{equation}\label{h-h_bar-sum}
\sum^{N-1}_{n=0} \nu_n =  \sum_{n=0}^{N-1} h^{'}_{t_n} Q_{t_n} \hat{Y}_{t_{n+1}} \tilde{\rho}_{t_n} \cdot \f{1}{2} \big( (V_{t_{n+1}} - V_{t_n})^2 - \Delta t \big) \stackrel{N \to \infty}{\longrightarrow} 0 \ \text{in } \,\L^1(\Omega, \tilde{\PP}_x).
\ee
In summary \eqref{dQ:term1} gives the estimate of the term ({\bf i}) in \eqref{Delta_n:step1} as 
\be \label{term-i}({\bf i})  =  \E_{0}^{x} \left[\int_{t_n}^{t_{n+1}} \hat{h}_s Q_s \d V_s \hat{Y}_{t_{n+1}} \right] +  \E_{0}^{x}[\nu_n] + \mathcal{O}\left((\Delta )^{\f{3}{2}}\right),\ee
with $\sum_{n=0}^{N-1} \E_{0}^{x}[\nu_n] \to 0$ in $\mathcal{L}^1(\Omega, \tilde{\PP}_x)$  as $N \to \infty$.

\bigskip
\noindent ({\bf ii})  It follows directly from the FBDSDEs system \eqref{FBDSDEs:interval} that term ({\bf ii}) in \eqref{Delta_n:step1} satisfies \begin{equation}\label{dQ:term2}
\begin{aligned}
({\bf ii})
 &=  \E_{0}^{x}\left[Q_{t_{n}}  \int_{t_n}^{t_{n+1}} \hat{Z}_s \d W_s - Q_{t_{n}} \int_{t_n}^{t_{n+1}} \left( \hat{h}_s \hat{Y}_s + \f{\tilde{\rho}_s}{\sigma_s} \hat{Z}_s \right) \d \ola{V}_s \right] \\[1.2ex]
 &= - \E_{0}^{x} \left[  Q_{t_{n}} \int_{t_n}^{t_{n+1}} \left( \hat{h}_s \hat{Y}_s  + \f{\tilde{\rho}_s}{\sigma_s} \hat{Z}_s \right) \d \ola{V}_s \right].
\end{aligned}
\end{equation}

\bigskip
\noindent ({\bf iii})  By splitting term  ({\bf iii})  in \eqref{Delta_n:step1}  and using the definition of $\eta_{n+1}$ in \eqref{def:eta} we obtain
\begin{equation}\label{eta:first}
\begin{aligned}
({\bf iii})  &= \E_{0}^{x}\left[ Q_{t_{n+1}}\hat{Y}^{'}_{t_{n+1}}  \eta_{n+1} \right] +  \f{1}{2}  \E_{0}^{x} \left [ Q_{t_{n+1}}  \hat{Y}^{''}_{t_{n+1}}\cdot ( \eta_{n+1} )^2 \right] + \E_0^x \left[Q_{t_{n+1}}  \xi_{n+1} \right] \\[1.2ex]
&= \underbrace{\E_{0}^{x}\left[ Q_{t_{n+1}} \hat{Y}^{'}_{t_{n+1}} \left( \int_{t_n}^{t_{n+1}} \rho_s \d W_s  -  \int_{t_n}^{t_{n+1}} \sigma_s \d W_s \right) \right]}_{({\bf iii - 1})} + \underbrace{\E_{0}^{x}\left[ Q_{t_n} \hat{Y}^{'}_{t_{n+1}}  \int_{t_n}^{t_{n+1}} \tilde{\rho}_s \d V_s \right]}_{({\bf iii - 2})}\\
&\quad  +  \underbrace{\E_{0}^{x}\left[ ( Q_{t_{n+1}} - Q_{t_n}) \hat{Y}^{'}_{t_{n+1}}  \int_{t_n}^{t_{n+1}} \tilde{\rho}_s \d V_s \right] -  \E_{0}^{x}\left[ Q_{t_{n+1}} \hat{Y}^{'}_{t_{n+1}}  \int_{t_n}^{t_{n+1}} \tilde{\rho}_s h_s \d s \right]}_{({\bf iii - 3})} \\
& \quad  +  \underbrace{\E_{0}^{x}\left[ Q_{t_{n+1}} \hat{Y}^{'}_{t_{n+1}}  R_X^{n+1} \right]}_{({\bf iii - 4})} +  \f{1}{2} \underbrace{\E_{0}^{x} \left [ Q_{t_{n+1}}  \hat{Y}^{''}_{t_{n+1}}\cdot ( \eta_{n+1} )^2 \right]}_{({\bf iii - 5})} + \E_0^x \left[Q_{t_{n+1}}  \xi_{n+1} \right].
\end{aligned}
\end{equation}
We next estimate terms ({\bf iii-1}) -- ({\bf iii-5}).  

\smallskip

Denote 
$$\nabla \hat{X}_{t} : = \f{\p \hat{X}_{t}^{t_n, x}}{\p x}|_{x = U_{t_n}}, \qquad \nabla \hat{Y}_{t} : = \f{\p Y_{t}^{t_n, x}}{\p x}|_{x = U_{t_n}}, \qquad \nabla \hat{Z}_{t} : = \f{\p Z_{t}^{t_n, x}}{\p x}|_{x = U_{t_n}}.$$
Then term ({\bf iii-1}) can be written as
\begin{equation}\label{eta:first:part1}
 ({\bf iii-1}) =  \E_{0}^{x}\left[ Q_{t_{n+1}} \left(\hat{Y}^{'}_{t_{n+1}} \nabla \hat{X}_{t_{n+1}} \right) \left( \int_{t_n}^{t_{n+1}} \rho_s \d W_s  -  \int_{t_n}^{t_{n+1}} \sigma_s \d W_s \right) \left(\nabla \hat{X}_{t_{n+1}}\right)^{-1} \right].
\end{equation}
By using the fact that $|( \nabla \hat{X}_{t_{n+1}})^{-1}| =1 + \mathcal{O}( \Delta t) $ and the
 following variational equation (see \cite{PP1994})
\begin{equation*}\label{nabla Y}
 \nabla \hat{Y}_{t} =  \hat{Y}^{'}_{t_{n+1}}  \nabla \hat{X}_{t_{n+1}} + \int_{t}^{t_{n+1}} \left( \hat{h}^{'}_s \hat{Y}_s  \nabla \hat{X}_s + \hat{h}_s \nabla \hat{Y}_s + \f{\tilde{\rho}_s}{\sigma_s} \nabla \hat{Z}_s \right) \d \ola{V}_s - \int_{t}^{t_{n+1}} \nabla \hat{Z}_s \d W_s,
\end{equation*}
we deduce that \eqref{eta:first:part1} becomes
\begin{equation*}
\begin{aligned}
({\bf iii-1}) &= \E_{0}^{x}\left[ Q_{t_{n+1}} \left( \hat{Y}^{'}_{t_{n+1}}  \nabla \hat{X}_{t_{n+1}} - \nabla \hat{Y}_{t_n}  \right) \left( \int_{t_n}^{t_{n+1}} \rho_s \d W_s  -  \int_{t_n}^{t_{n+1}} \sigma_s \d W_s \right)  \right] \\
& \quad + \E_{0}^{x}\left[ Q_{t_{n+1}} \nabla \hat{Y}_{t_n} \left( \int_{t_n}^{t_{n+1}} \rho_s \d W_s  -  \int_{t_n}^{t_{n+1}} \sigma_s \d W_s \right)\right] + \mathcal{O}((\Delta t)^{\f{3}{2}})\\[1.2ex]
& = \E_{0}^{x}\left[ \left(Q_{t_{n}}  \int_{t_n}^{t_{n+1}} \nabla \hat{Z}_s \d W_s + \lambda_{t_n}\right)  \cdot \left( \int_{t_n}^{t_{n+1}} \rho_s \d W_s  -  \int_{t_n}^{t_{n+1}} \sigma_s \d W_s \right) \right] + \mathcal{O}((\Delta t)^{\f{3}{2}}),
\end{aligned}
\end{equation*}
where $\lambda_{t_n} = - Q_{t_n}  \int_{t_n}^{t_{n+1}} \big( \hat{h}^{'}_s \bar{Y}_s  \nabla \hat{X}_s + \hat{h}_s \nabla \hat{Y}_s + \f{\tilde{\rho}_s}{\sigma_s} \nabla \hat{Z}_s \big)d \ola{V}_s  + Q_{t_{n+1} }  \nabla \hat{Y}_{t_n}$ is independent of $\int_{t_n}^{t_{n+1}} \rho_s \d W_s  -  \int_{t_n}^{t_{n+1}} \sigma_s \d W_s$ and hence gives
$$ \E_{0}^{x}\left[ \lambda_{t_n} \left( \int_{t_n}^{t_{n+1}} \rho_s \d W_s  -  \int_{t_n}^{t_{n+1}} \sigma_s \d W_s \right) \right] = 0. $$
As a consequence
\begin{equation}\label{estiii-1}
\begin{aligned}
({\bf iii-1}) & =   \E_{0}^{x}\left[ Q_{t_{n}} \int_{t_n}^{t_{n+1}} \nabla \hat{Z}_s \d W_s \cdot  \left( \int_{t_n}^{t_{n+1}} \rho_s \d W_s  -  \int_{t_n}^{t_{n+1}} \sigma_s \d W_s \right) \right]
+ \mathcal{O}\left((\Delta t)^{\f{3}{2}}\right) \\[1.2ex]
&= \E_{0}^{x}\left[ Q_{t_{n}} \nabla \hat{Z}_{t_{n+1}}  \cdot \int_{t_n}^{t_{n+1}}  (\rho_s - \sigma_s )\d s \right] + \mathcal{O}\left((\Delta t)^{\f{3}{2}}\right).
\end{aligned}
\end{equation}

\smallskip
 Let $C$ represent a generic constant while the context is clear.    By the definition of $R_X^{n+1}$, it is straightforward to verify that 
 \be\label{estiii-4} ({\bf iii-4}) = \E_{0}^{x}\left[ Q_{t_{n+1}} \hat{Y}^{'}_{t_{n+1}}  R_X^{n+1} \right] \leq C (\Delta t)^{\f{3}{2}}. \ee 
 
 \smallskip
 
  Applying It\^o formula to $Q_{t_n}$ in term ({\bf iii-3}) we obtain
\be\label{eta:first:est2}
\begin{aligned}
({\bf iii-3}) & = \E_{0}^{x}\left[ ( Q_{t_{n+1}} - Q_{t_n}) \hat{Y}^{'}_{t_{n+1}}  \int_{t_n}^{t_{n+1}} \tilde{\rho}_s \d V_s \right] -  \E_{0}^{x}\left[ Q_{t_{n+1}} \hat{Y}^{'}_{t_{n+1}}  \int_{t_n}^{t_{n+1}} \tilde{\rho}_s h_s \d s \right]\\[1.2ex]
& =  \E_{0}^{x}\left[ \int_{t_n}^{t_{n+1}} h_s Q_s \d V_s  \int_{t_n}^{t_{n+1}} \tilde{\rho}_s \d V_s \hat{Y}^{'}_{t_{n+1}} \right]  - \E_{0}^{x}\left[ Q_{t_{n+1}} \hat{Y}^{'}_{t_{n+1}}  \int_{t_n}^{t_{n+1}} \tilde{\rho}_s h_s \d s \right] \\[1.2ex]
& \leq  \left| \E_{0}^{x}\left[ \int_{t_n}^{t_{n+1}} h_s ( Q_s - Q_{t_{n+1}} ) \d V_s  \int_{t_n}^{t_{n+1}} \tilde{\rho}_s \d V_s \hat{Y}^{'}_{t_{n+1}} \right] \right| \\
& \qquad + \left| \E_{0}^{x}\left[ Q_{t_{n+1}} \hat{Y}^{'}_{t_{n+1}} \left( \int_{t_n}^{t_{n+1}} h_s  \d V_s  \int_{t_n}^{t_{n+1}} \tilde{\rho}_s \d V_s  -   \int_{t_n}^{t_{n+1}} \tilde{\rho}_s h_s \d s \right) \right] \right | \\
& \leq  C (\Delta t)^{\f{3}{2}}.
\end{aligned}
\end{equation}

\smallskip


By using the definition of $\eta_{n+1}$ in \eqref{def:eta}, we deduce that
\begin{equation*}
\begin{aligned}
({\bf iii-5}) &=  \f{1}{2} \E_0^x\left[ Q_{t_n}  \hat{Y}^{''}_{t_{n+1}}  \int_{t_{n}}^{t_{n+1}} \left(  \rho^2_s + \tilde{\rho}_s^2 + \sigma_s^2 - 2 \rho_s \sigma_s \right) \d s \right] + \mathcal{O}\left((\Delta t)^{\f{3}{2}}\right) \\
&=  \E_0^x \left[Q_{t_n}  \hat{Y}^{''}_{t_{n+1}}  \int_{t_{n}}^{t_{n+1}} \left(  \sigma_s^2 - \rho_s \sigma_s \right) \d s \right] + \mathcal{O}\left((\Delta t)^{\f{3}{2}}\right). 
\end{aligned}
\end{equation*}
As a simple corollary of the assertion \eqref{Z=dY}, we have $\hat{Y}^{''}_{t_{n+1}} \sigma_{t_{n+1}}= \nabla \hat{Z}_{t_{n+1}} + \mathcal{O}(\Delta t)$ and thus 
\be \label{eta:second}
({\bf iii-5}) = \E_0^x \left[Q_{t_n}  \nabla \hat{Z}_{t_{n+1}}  \int_{t_{n}}^{t_{n+1}} \left(  \sigma_s - \rho_s \right) \d s \right] +\mathcal{O}\left((\Delta t)^{\f{3}{2}}\right).   
\ee

\smallskip
It then remains to estimate term ({\bf iii-2}).     Notice that due to equations  \eqref{Z=dY} and \eqref{nabla Y} we have ${\hat{Z}_{s}}/{\sigma_{s}} = \nabla \hat{Y}_{s} (\nabla \hat{X}_s)^{-1}$.   Hence for any $s \in [t_n, t_{n+1}]$ it holds
$$
\begin{aligned}
\hat{Y}^{'}_{t_{n+1}}  - \f{\hat{Z}_{s}}{\sigma_{s}} & = \hat{Y}^{'}_{t_{n+1}}  - \nabla \hat{Y}_{s} (\nabla \hat{X}_s)^{-1} \\
& = - \int_{s}^{t_{n+1}} \left( \hat{h}^{'}_r \hat{Y}_r \nabla \hat{X}_r + \hat{h}_r \nabla \hat{Y}_r + \f{\tilde{\rho}_r}{\sigma_r} \nabla \hat{Z}_r \right) \d \ola{V}_r - \int_{s}^{t_{n+1}} \nabla \hat{Z}_r \d W_r + \mathcal{O}(\Delta t), 
\end{aligned} 
$$
and therefore
$$
\begin{aligned}
- Q_{t_{n}} \int_{t_n}^{t_{n+1}} \f{\tilde{\rho_s}}{\sigma_s} \hat{Z}_s  \d \ola{V}_s &= - Q_{t_{n} }  \hat{Y}^{'}_{t_{n+1}}  \int_{t_n}^{t_{n+1}} \tilde{\rho_s} \d \ola{V}_s \\
& \quad - Q_{t_{n} }  \int_{t_n}^{t_{n+1}} \tilde{\rho_s} \int_{s}^{t_{n+1}} \left( \hat{h}^{'}_r \hat{Y}_r \nabla \hat{X}_r + \hat{h}_r \nabla \hat{Y}_r + \f{\tilde{\rho}_r}{\sigma_r} \nabla \hat{Z}_r \right) \d \ola{V}_r \d \ola{V}_s \\
& \quad - Q_{t_{n} }  \int_{t_n}^{t_{n+1}} \tilde{\rho_s} \int_{s}^{t_{n+1}} \nabla \hat{Z}_r \d W_r \d \ola{V}_s + \mathcal{O}\left((\Delta t)^{\f{3}{2}}\right).
\end{aligned}
$$
Since $W$ and $V$ are two independent Brownian motions, 
\begin{equation}\label{dZ-dW-dV}
\begin{aligned}
 & \E_{0}^{x} \left[ Q_{t_{n} }  \int_{t_n}^{t_{n+1}} \tilde{\rho_s} \int_{s}^{t_{n+1}} \nabla \hat{Z}_r \d W_r \d \ola{V}_s \right] \\[1.2ex]
& =  \E_{0}^{x} \left[ Q_{t_{n} }  \int_{t_n}^{t_{n+1}} \tilde{\rho_s} \int_{s}^{t_{n+1}} ( \nabla \hat{Z}_r - \nabla \hat{Z}_{t_{n}}) \d W_r \d \ola{V}_s \right]  \leq C(\Delta t)^{\f{3}{2}}.
\end{aligned}
\end{equation}
As a result, 
\begin{equation}\label{estiii-2}
 \begin{aligned}
 ({\bf iii-2}) & = \E_{0}^{x} \left[ Q_{t_{n}} \int_{t_n}^{t_{n+1}} \f{\tilde{\rho}_s}{\sigma_s} \hat{Z}_s  \d \ola{V}_s \right]\\ 
  &\quad - \E_{0}^{x} \left[ Q_{t_{n} }  \int_{t_n}^{t_{n+1}} \tilde{\rho}_s \int_{s}^{t_{n+1}} \left( \hat{h}^{'}_r \hat{Y}_r \nabla \hat{X}_r + \hat{h}_r \nabla \hat{Y}_r + \f{\tilde{\rho}_r}{\sigma_r} \nabla \hat{Z}_r \right) \d \ola{V}_r \d \ola{V}_s \right] + \mathcal{O}\left( (\Delta t)^{\f{3}{2}}\right) \\[1.2ex]
  & =  \E_{0}^{x} \left[ Q_{t_{n}} \int_{t_n}^{t_{n+1}} \f{\tilde{\rho}_s}{\sigma_s} \hat{Z}_s  \d \ola{V}_s \right]  -  \E_{0}^{x}\left[ \tilde{\lambda}_{t_n}  \int_{t_n}^{t_{n+1}}  \int_{s}^{t_{n+1}} \d \ola{V}_r \d \ola{V}_s \right] + \mathcal{O}\left( (\Delta t)^{\f{3}{2}}\right),
  \end{aligned}
 \end{equation}
where $ \tilde{\lambda}_{t_n} = Q_{t_{n} } \tilde{\rho}_{t_{n+1}}\left( \hat{h}^{'}_{t_n} \hat{Y}_{t_{n+1}} \nabla \hat{X}_{t_n} + \hat{h}_{t_n} \nabla \hat{Y}_{t_{n+1}} + \f{ \tilde{\rho}_{t_{n+1}} }{\sigma_{t_{n+1}}} \nabla \hat{Z}_{t_{n+1}} \right)$
 is independent of $ \int_{t_n}^{t_{n+1}}  \int_{s}^{t_{n+1}} d \ola{V}_r d \ola{V}_s = \f{1}{2} \big( (V_{t_{n+1}} - V_{t_n})^2 - \Delta t \big)$.   By an argument similar to \eqref{h-h_bar-sum}, we obtain
 \be\label{estiii-2.1}
 \sum_{n=0}^{N-1} \E_{0}^{x} \left[ \tilde{\lambda}_{t_n}  \int_{t_n}^{t_{n+1}}  \int_{s}^{t_{n+1}} \d \ola{V}_r \d \ola{V}_s \right]  \stackrel{N \to \infty}{\longrightarrow} 0 \ \text{in } \mathcal{L}^1(\Omega, \tilde{\PP}_x). \ee 

\smallskip
Collecting estimates \eqref{estiii-1}, \eqref{estiii-4}, \eqref{eta:first:est2}, \eqref{eta:second} and \eqref{estiii-2} into \eqref{eta:first}; then inserting \eqref{eta:first}, \eqref{term-i} and \eqref{dQ:term2} into 
\eqref{Delta_n:step1} we finally obtain 
\begin{equation}\label{a-b-c}
\begin{aligned}
\Delta_n &=   \E_{0}^{x}\left[ \int_{t_n}^{t_{n+1}} \hat{h}_s Q_s \d V_s \hat{Y}_{t_{n+1}} - Q_{t_{n}} \int_{t_n}^{t_{n+1}} \hat{h}_s \hat{Y}_s \d \ola{V}_s  + \nu_n \right]+\mathcal{O}\left( (\Delta t)^{\f{3}{2}}\right) \\[1.2ex]& =   \E_{0}^{x}[\alpha_n] + \E_{0}^{x}[\beta_n] + \E_{0}^{x}[\gamma_n] + \E_{0}^{x}[\nu_n] + \mathcal{O}\left( (\Delta t)^{\f{3}{2}}\right) ,
\end{aligned}
\end{equation}
where $\{\nu_n\}$ is defined as in \eqref{def:nu} satisfying \eqref{h-h_bar-sum},
and
\begin{eqnarray*}
\alpha_n &: =&  \ \int_{t_n}^{t_{n+1}} \left( Q_s \hat{h}_s - Q_{t_n} \hat{h}_{t_{n}}\right)\hat{Y}_{t_{n+1}} \d V_s, \\
\beta_n &:=&  \int_{t_n}^{t_{n+1}}  Q_{t_n}\left(  \hat{h}_{t_{n+1}} \hat{Y}_{t_{n+1}} -  \hat{h}_s \hat{Y}_s  \right) \d \ola{V}_s, \\
\gamma_n &:=& Q_{t_n} \hat{Y}_{t_{n+1}} \left(  \hat{h}_{t_{n}} -  \hat{h}_{t_{n+1}} \right) \cdot ( V_{t_{n+1}} - V_{t_n} ). \end{eqnarray*}
The last steps are to show that $\sum_{n=0}^{N-1} \E_{0}^{x}[\alpha_n] \to 0$, $\sum_{n=0}^{N-1} \E_{0}^{x}[\beta_n] \to 0$, and $\sum_{n=0}^{N-1} \E_{0}^{x}[\gamma_n] \to 0$ in $\mathcal{L}^1(\Omega, \tilde{\PP}_x)$  as $N \to \infty$.  

\smallskip 
First write $\alpha_n = \alpha_n^{(1)}+ \alpha_n^{(2)} + \alpha_n^{(3)}$ with
\begin{eqnarray*}
\alpha_n^{(1)} &:=&   \int_{t_n}^{t_{n+1}} (Q_s - Q_{t_n} )  \hat{h}_s \d V_s \cdot \hat{Y}_{t_{n+1}},  \\
 \alpha_n^{(2)} &:=&  \int_{t_n}^{t_{n+1}} Q_{t_n} ( \hat{h}_s -  \hat{h}_{t_n}) \d V_s \cdot \hat{Y}_{t_{n}},\\
  \alpha_n^{(3)} &:=&  \int_{t_n}^{t_{n+1}} Q_{t_n} ( \hat{h}_s -  \hat{h}_{t_n}) \d V_s \cdot \left( \hat{Y}_{t_{n+1}} - \hat{Y}_{t_{n}} \right).  
\end{eqnarray*}
Denote by $\tilde{\E}_{x} $ the expectation with respect to $\tilde{\PP}_{x}$, where $\tilde{\PP}_x:= \tilde{\PP}(\cdot | U_0 = x)$ is the induced probability measure. 
Notice that $\hat{Y}_{t_n} = Y_{t_n}(U_{t_n})$ due to $\hat{X}_{t_n} = U_{t_n}$ given in \eqref{FBDSDEs:interval}, and that $h_t$ is a bounded function, we apply It\^o's formula to $(Q_s - Q_{t_n} )$ in $\alpha_n^{(1)}$ to get
$$
\begin{aligned}
\tilde{\E}_x \left[ \left|\E_{0}^{x}\big[ \sum_{n = 0}^{N-1} \alpha_n^{(1)} \big] \right| \right]  = & \ \tilde{\E}_x\left[ \left|  \E_{0}^{x}\big[ \sum_{n = 0}^{N-1} \int_{t_n}^{t_{n+1}} \int_{t_n}^s  \hat{h}_r Q_r \d V_r \hat{h}_s \d V_s \cdot \hat{Y}_{t_{n+1}} \big] \right| \right]  \\
\leq & \ \tilde{\E}_x\left[ \left| \sum_{n = 0}^{N-1} \int_{t_n}^{t_{n+1}} \int_{t_n}^s  \hat{h}_r ( Q_r - Q_{t_n} )  \d V_r  \hat{h}_s \d V_s \cdot \hat{Y}_{t_{n+1}}  \right|  \right] \\ 
& \quad + \tilde{\E}_x\left[ \left| \sum_{n = 0}^{N-1} \int_{t_n}^{t_{n+1}} \int_{t_n}^s  \hat{h}_{r} Q_{t_n} \d V_r  \hat{h}_s \d V_s \cdot \bar{Y}_{t_{n+1}} \right|  \right]  \\
\leq & \ C \sum_{n = 0}^{N-1} (\Delta t)^{\f{3}{2}}  + C  \tilde{\E}_x\left[ \left| \sum_{n = 0}^{N-1} Q_{t_n} \bar{Y}_{t_{n+1}} \cdot   \f{1}{2} \left( (V_{t_{n+1}} - V_{t_n})^2 - \Delta t\right)  \right|  \right]  .
\end{aligned}
$$ 
and from the  fact that   
$$\sum_{n = 0}^{N-1} Q_{t_n} \bar{Y}_{t_{n+1}}  \cdot \f{1}{2} \big( (V_{t_{n+1}} - V_{t_n})^2 - \Delta t\big)  \rightarrow 0 \ \text{in} \ \mathcal{L}^1(\Omega, \tilde{P}_x) $$ 
we have 
\begin{equation}\label{converge:alpha_n^1}
\E_{0}^{x}\big[ \sum_{n=0}^{N-1} \alpha_n^{(1)} \big] \rightarrow 0 \ \text{in} \ \mathcal{L}^1(\Omega, \tilde{P}_x).
\end{equation}
For $ \alpha_n^{(2)}$, we apply It\^o's formula to $h_s$ to get
\begin{equation*}
\begin{aligned}
\E_{0}^x[  \alpha_n^{(2)}  ]  = & \E_{0}^x\big[    \int_{t_n}^{t_{n+1}} Q_{t_n} \cdot \left(  \int_{t_n}^{s}  [ \hat{b}_r \cdot \hat{h}^{\prime}_r + \f{ (\sigma_r)^2 }{2} \cdot \hat{h}^{\prime \prime}_r ] dr + \int_{t_n}^{s} \sigma_r \cdot \hat{h}^{\prime}_r d W_r  \right) \cdot \hat{Y}_{t_{n}} dV_s   \big]  \\
= & \E_{0}^x\big[ \int_{t_n}^{t_{n+1}} Q_{t_n} \hat{Y}_{t_{n}}  \cdot \int_{t_n}^{s}   [ \hat{b}_r \cdot \hat{h}^{\prime}_r + \f{ (\sigma_r)^2 }{2} \hat{h}^{\prime \prime}_r ]  dr dV_s   \big].  \\
\end{aligned}
\end{equation*}
Since $b_t$, $\sigma_t$, $h^{\prime}$ and $h^{\prime \prime}$ are all bounded,
we have
\begin{equation*}
\begin{aligned}
\tilde{E}_x[ | \E_{0}^x[  \alpha_n^{(2)}  ]  | ]  \leq & C \Delta t^{3/2}.
\end{aligned}
\end{equation*}
Moreover, it follows from H$\ddot{{\rm o}}$lder's inequality, Lemma \ref{Y:regularity} and Lemma \ref{h:regularity} that
\begin{equation*}
\begin{aligned}
 \tilde{E}_x[ | \E_{0}^x[ \alpha_n^{(3)} ] | ]  = &   \tilde{E}_x\Big[ \big| \E_{0}^x\big[\int_{t_n}^{t_{n+1}} Q_{t_n} (\hat{h}_s - \hat{h}_{t_n}) d V_s \cdot \big( \hat{Y}_{t_{n+1}} - \hat{Y}_{t_{n}} \big) \big]\big|  \Big]  \\
 \leq &  \Big(\tilde{E}_x\big[ \big( \int_{t_n}^{t_{n+1}} Q_{t_n} ( \hat{h}_s -  \hat{h}_{t_n}) d V_s \big) ^2 \big] \Big)^{\f{1}{2}}  \cdot \Big(\tilde{E}_x\big[  \big( \hat{Y}_{t_{n+1}} - \hat{Y}_{t_{n}} \big)^2 \big]\Big)^{\f{1}{2}}  \\
 \leq & C \Delta t^{3/2}.
 \end{aligned}
\end{equation*}
Hence,
\begin{equation}\label{converge:alpha_n^2}
 = \sum_{n=0}^{N-1} \E_{0}^x[ \alpha_n^{(2)}]  \rightarrow 0 \ \text{in} \ \mathcal{L}^1(\Omega, \tilde{P}_x),
\end{equation}
and
\begin{equation}\label{converge:alpha_n^3}
\sum_{n=0}^{N-1}  \E_{0}^x[\alpha_n^{(3)} ]\rightarrow 0 \ \text{in} \ \mathcal{L}^1(\Omega, \tilde{P}_x).
\end{equation}
Then, from \eqref{converge:alpha_n^1}, \eqref{converge:alpha_n^2} and \eqref{converge:alpha_n^3}, we get
\begin{equation}\label{converge:alpha_n}
\sum_{n=0}^{N-1}  \E_{0}^x[\alpha_n]  \rightarrow 0 \ \text{in} \ \mathcal{L}^1(\Omega, \tilde{P}_x).
\end{equation}

For the term $\beta_n$  in \eqref{a-b-c}, we have 
$$
\begin{aligned}
\beta_n = & \int_{t_n}^{t_{n+1}}  Q_{t_n}\Big( ( \hat{h}_{t_{n+1}} -  \hat{h}_s )  \hat{Y}_{t_{n+1}} +   \hat{h}_s (  \hat{Y}_{t_{n+1}} -  \hat{Y}_s ) \Big) d\ola{V}_s \\
= & \int_{t_n}^{t_{n+1}}  Q_{t_n} (  \hat{h}_{t_{n+1}} -   \hat{h}_s )  \hat{Y}_{t_{n}} d\ola{V}_s+ \int_{t_n}^{t_{n+1}}  Q_{t_n}  \hat{h}_s (  \hat{Y}_{t_{n+1}} -  \hat{Y}_s )  d\ola{V}_s \\
& \quad + \int_{t_n}^{t_{n+1}}  Q_{t_n} (  \hat{h}_{t_{n+1}} -  \hat{h}_s ) d\ola{V}_s \cdot \big(  \hat{Y}_{t_{n+1}} -  \hat{Y}_{t_{n}}  \big)  \\
= & \beta_n^{(1)} + \beta_n^{(2)} + \beta_n^{(3)}
\end{aligned}
$$
with
$$\beta_n^{(1)} = \int_{t_n}^{t_{n+1}}  Q_{t_n} (  \hat{h}_{t_{n+1}} -   \hat{h}_s )  \hat{Y}_{t_{n}} d\ola{V}_s ,$$
$$\beta_n^{(2)} = \int_{t_n}^{t_{n+1}}  Q_{t_n}   \hat{h}_s (  \hat{Y}_{t_{n+1}} -  \hat{Y}_s )  d\ola{V}_s $$
and
$$ \beta_n^{(3)} =  \int_{t_n}^{t_{n+1}}  Q_{t_n} (  \hat{h}_{t_{n+1}} -   \hat{h}_s ) d\ola{V}_s  \cdot \big(  \hat{Y}_{t_{n+1}} -  \hat{Y}_{t_{n}}  \big) . $$
Following the similar the approaches to $\alpha_n^{(2)}$ and $\alpha_n^{(3)}$, we have
\begin{equation*}
\begin{aligned}
\tilde{E}_x[ |\E_{0}^x[\beta_n^{(1)} ] | ] \leq C \Delta t^{3/2}
\end{aligned}
\end{equation*}
and
\begin{equation*}
\begin{aligned}
\tilde{E}_x[ | \E_{0}^x[\beta_n^{(3)}] | ] \leq C \Delta t^{3/2}.
\end{aligned}
\end{equation*}
Hence,
\begin{equation}\label{converge:beta_n^1}
\sum_{n=0}^{N-1} \E_{0}^x[\beta_n^{(1)} ] \rightarrow 0 \ \text{in} \ \mathcal{L}^1(\Omega, \tilde{P}_x),
\end{equation}
and
\begin{equation}\label{converge:beta_n^3}
\sum_{n=0}^{N-1} \E_{0}^x[ \beta_n^{(3)} ] \rightarrow 0 \ \text{in} \ \mathcal{L}^1(\Omega, \tilde{P}_x).
\end{equation}
From the BDSDE in \eqref{Backward:FBDSDEs}, Lemma \ref{Y:regularity}, Lemma \ref{Y:bounded} and estimate \eqref{dZ-dW-dV}, we get
\begin{equation*}
\begin{aligned}
\sum_{n=0}^{N-1}\E_0^x[ \beta_n^{(2)} ] = & \sum_{n=0}^{N-1}\E_0^x\Big[ \int_{t_n}^{t_{n+1}}  Q_{t_n}  \hat{h}_s \Big( \int_{s}^{t_{n+1}} \hat{Z}_r dW_r - \int_{s}^{t_{n+1}} \big( \hat{h}_r \hat{Y}_r + \f{\tilde{\rho}_r }{\sigma_r} \hat{Z}_r  \big) d\ola{V}_r \Big)  d\ola{V}_s   \Big] \\
= & \sum_{n=0}^{N-1} \E_0^x\Big[ \int_{t_n}^{t_{n+1}}  Q_{t_n}  \hat{h}_s \big( - \int_{s}^{t_{n+1}}  \big( \hat{h}_r \hat{Y}_{t_{n+1}} + \f{\tilde{\rho}_r }{\sigma_r} \hat{Z}_{t_{n+1}} \big) d\ola{V}_r \big)  d\ola{V}_s  \Big]  + O((\Delta t)^{\f{3}{2}}) \\ 
 & \quad + \sum_{n=0}^{N-1}\E_0^x\Big[ \int_{t_n}^{t_{n+1}}  Q_{t_n}  \hat{h}_s \big( \int_{s}^{t_{n+1}}  \hat{h}_r ( \hat{Y}_{t_{n+1}}  - \hat{Y}_r  ) + \f{\tilde{\rho}_r }{\sigma_r} (\hat{Z}_{t_{n+1}} - \hat{Z}_r )  d\ola{V}_r \big)  d\ola{V}_s  \Big]  \\ 
\end{aligned}
\end{equation*}
Next, we take conditional expectation $\tilde{E}_x$ to the absolute value of the above equation. Since
$$\tilde{E}_x\big[ ( \hat{Y}_{t_{n+1}} - \hat{Y}_r  )^2 \big] \leq C(\Delta t)^{3}, \ \tilde{E}_x\big[ ( \hat{Z}_{t_{n+1}} - \hat{Z}_r  )^2 \big] \leq C(\Delta t)^{3}$$
and
$$
\begin{aligned}
& \tilde{E}_x\Big[ \big|  \E_0^x \big[\sum_{n=0}^{N-1}  \int_{t_n}^{t_{n+1}}  Q_{t_n}  \hat{h}_s \big( - \int_{s}^{t_{n+1}}  \hat{h}_r \hat{Y}_{t_{n+1}} d\ola{V}_r \big)  d\ola{V}_s  \big] \big| \Big] \\
\leq &  C \tilde{E}_x\Big[ \big| \sum_{n=0}^{N-1} Q_{t_n} \hat{Y}_{t_{n+1}}  \cdot \f{1}{2} \big( (V_{t_{n+1}} - V_{t_n})^2 - \Delta t\big)   \big| \Big],
\end{aligned} 
 $$
it follows from the fact 
$$\sum_{n = 0}^{N-1} Q_{t_n} \hat{Y}_{t_{n+1}}  \cdot \f{1}{2} \big( (V_{t_{n+1}} - V_{t_n})^2 - \Delta t\big)  \rightarrow 0 \ \text{in} \ \mathcal{L}^1(\Omega, \tilde{P}_x) $$ 
that
\begin{equation}\label{converge:beta_n^2}
\E_0^x[ \sum_{n=0}^{N-1} \beta_n^{(2)} ]  \rightarrow 0 \ \text{in} \ \mathcal{L}^1(\Omega, \tilde{P}_x).
\end{equation}
Hence,
\begin{equation}\label{converge:beta_n}
\E_0^x[\sum_{n=0}^{N-1} \beta_n ] \rightarrow 0 \ \text{in} \ L^1(\Omega, \tilde{P}_x).
\end{equation}

For $\gamma_n$, applying It\^o formula to $h_t$, it's easy to verify that
\begin{equation}\label{estimate:E[h]}
| \E_0^x[ \hat{h}_{t_{n+1}} - \hat{h}_{t_n}] | \leq C \Delta t.
\end{equation}
Since $ \hat{h}_{t_n} - \hat{h}_{t_{n+1}}$ is independent from $ Q_{t_n} \hat{Y}_{t_n}( V_{t_{n+1}} - V_{t_n} ) $,
\begin{equation*}
\begin{aligned}
\E_0^x[  \gamma_n ]  =  & \ \E_0^x\big[Q_{t_n} \hat{Y}_{t_{n}} \big(  \hat{h}_{t_{n}} -  \hat{h}_{t_{n+1}} \big) \cdot ( V_{t_{n+1}} - V_{t_n} ) \big] \\
& + \E_0^x\big[Q_{t_n} \big( \hat{Y}_{t_{n+1}} - \hat{Y}_{t_{n}} \big) \cdot \big(  \hat{h}_{t_{n}} -  \hat{h}_{t_{n+1}} \big) \cdot ( V_{t_{n+1}} - V_{t_n} ) \big] \\
= & \ \E_0^x\big[ \hat{h}_{t_{n}} -  \hat{h}_{t_{n+1}} \big] \cdot \E_0^x\big[Q_{t_n} \hat{Y}_{t_{n}} \cdot ( V_{t_{n+1}} - V_{t_n} ) \big] \\
& + \E_0^x\big[Q_{t_n} \big( \hat{Y}_{t_{n+1}} - \hat{Y}_{t_{n}} \big) \cdot \big(  \hat{h}_{t_{n}} - \hat{h}_{t_{n+1}} \big) \cdot ( V_{t_{n+1}} - V_{t_n} ) \big] .
\end{aligned}
\end{equation*}
Then, from estimate \eqref{estimate:E[h]}, lemma \ref{Y:regularity}, lemma \ref{h:regularity}, we get
\begin{equation}\label{convergence:gamma_n}
\begin{aligned}
\tilde{E}_x \big[ | \E_0^x[  \gamma_n ]  | \big] \leq & \ C \Delta t \cdot \tilde{E}_x \big[ \big| Q_{t_n} \hat{Y}_{t_{n}} \cdot ( V_{t_{n+1}} - V_{t_n} ) \big| \big] \\
& + \tilde{E}_x\big[ | Q_{t_n} \big( \hat{Y}_{t_{n+1}} - \hat{Y}_{t_{n}} \big) \cdot \big(  \hat{h}_{t_{n}} -  \hat{h}_{t_{n+1}} \big) \cdot ( V_{t_{n+1}} - V_{t_n} ) | \big] \\
\leq & C (\Delta t)^{\f{3}{2}}
\end{aligned}
\end{equation}
and therefore
\begin{equation}\label{converge:gamma_n}
\E_0^x[ \sum_{n=0}^{N-1} \gamma_n ] \rightarrow \ \text{in} \ L^1(\Omega, \tilde{P}_x).
\end{equation}

Finally with convergence results in \eqref{converge:alpha_n}, \eqref{converge:beta_n} and \eqref{converge:gamma_n}, we have
$$  \sum_{n=0}^{N-1} \Delta_n \rightarrow 0 \ \text{in} \ L^1(\Omega, \tilde{P}_x)  $$
as required.  \hfill$\Box$

\bigskip
\subsection{Adjoint FBDSDEs}

In this subsection, we consider the following FBDSDEs system, in which the ``forward SDE''   \eqref{FBDSDEs:eq1} goes backward and the ``Backward SDE'' \eqref{FBDSDEs:eq2} goes forward
\begin{equation}\label{Forward:FBDSDEs}
\left\{
\begin{aligned}
\d \ola{X}_t  =&  b_t (\ola{X}_t) \d t - \sigma_t \d \ola{W}_t, \quad 0 \leq t \leq \tau & \text{(SDE)} \\
\d \ora{Y}_t = & - b_t^{\prime}(\ola{X}_t) \ora{Y}_t  \d t -  \ora{Z}_t^{T, x} \d \ola{W}_t + \left( h(\ola{X}_t) \ora{Y}_t  - \f{\tilde{\rho}_t}{\sigma_t} \ora{Z}_t  \right) \d V_t, \quad 0 \leq t \leq \tau  & \text{(BDSDE)} \\
\ola{X}_{\tau} =& x, \qquad 
\ora{Y}_0 = p_0(\ola{X}_0), 
\end{aligned}\right.
\end{equation}
where $0\leq \tau \leq T$, $\int_{t}^{T} \cdot \d \ola{W}_s$ is a backward It\^o Integral and $ \int_{t}^{T} \cdot \d V_s$ is a standard forward It\^o integral.   Write the solution to \eqref{Forward:FBDSDEs} as $( \ola{X}_t^{T, x}, \ora{Y}_t^{T, x}, \ora{Z}_t^{T, x})$.   Then
by inverting the time index in the standard FBDSDEs system,  $\ola{X}_t^{T, x}$ is a $\mathcal{F}^{W}_{t, T}$ adaptive stochastic process and the solution $\left( \ora{Y}_t^{T, x}, \ora{Z}_t^{T, x}\right)$ of the BDSDE in \eqref{Forward:FBDSDEs} is adaptive to $ \mathcal{F}^{W}_{t, T} \vee \mathcal{F}^{V}_t $.    Similar to the notation used in Section 2.1, we denote $\ora{Y}_{t}(x) : = \ora{Y}_{t}^{t, x}$ and $\ora{Z}_{t}(x) : = \ora{Z}_{t}^{t, x}$.  In addition, for any non-negative integer $m$ and function $\eta_t(x)$ we write $\eta_t^{(m)}: = \f{\p^{m} \eta}{\p x^{m}}$.

We need the following regularity properties for $b$ and $\sigma$. 
\begin{ass}\label{b-sigma:regularity}
For $0 \leq s \leq t \leq T$,  functions $b$ and $\sigma$ satisfy
$$\left|b_t(x) - b_s(x)\right| + \left|b^{\prime}_t(x) - b^{\prime}_s(x)\right| \leq C | t - s |, \qquad
\left|\sigma_t - \sigma_s \right| \leq C |t - s|,$$
where $C$ is a given positive constant independent of $b$, $\sigma$, $s$ and $t$.
\end{ass}

\smallskip

Lemma \ref{Y:bounded-support} can be proved by using  repeatedly the variational form of BDSDEs \cite{PP1994}.   
\begin{lem}\label{Y:bounded-support}
Assume that $b \in \C_{l, b}^{4}$, $\phi \in \C_{l, b}^{3}$, $h \in \C_{l, b}^{3}$ and every derivative of $b$, $\phi$ and $h$ has bounded support in $\mathbb{R}$.  Then for each $m_1 = 0, 1, 2$ and $m_2 = 0, 1, 2, 3$,   $Y_t^{(m_1)}$, $\ora{Y}^{(m_2)}_t$ have bounded support and satisfy
\begin{equation}\label{Y:integral-bounded}
\int_{\mathbb{R}}  \tilde{\E}\left[ \sup_{0 \leq t \leq T} \left| Y_t^{(m_1)} \right|^2 \right] \d x  < \infty \  \text{and}  \ \int_{\mathbb{R}} \tilde{\E}\left[ \sup_{0 \leq t \leq T}\left| \ora{Y}_t^{(m_2)} \right|^2 \right] \d x  < \infty.
\end{equation}
\end{lem}

\smallskip
Denote by $\la \cdot , \cdot  \ra$ the standard inner product in $\L^2$.  The following theorem shows that $\ora{Y}_t$ is the adjoint stochastic process of $Y_t$ defined in the FBDSDEs system \eqref{Backward:FBDSDEs}.

\begin{thm}\label{thm:Adjoint-Operator}
Assume that, in addition to Assumption \ref{b-sigma:regularity}holds, $\sigma$ is uniformly bounded, $b \in \C_{l, b}^{4}$, $\phi \in \C_{l, b}^{3}$, $h \in \C_{l, b}^{3}$ and each derivative of $b$, $\phi$ and $h$ has bounded support in $\mathbb{R} $.  Then the process $R_t := \la Y_{t}, \ora{Y}_{t} \ra$, $t \in [0, T]$ is a constant  for almost all trajectories.
\end{thm}

\smallskip

\noindent {\it Proof.}    According to \cite{PP1994},  $R_t$ has a.s. continuous paths, it suffices to show that $\forall s, t \in [0, T]$, $R_s = R_t$ a.s.. 

\smallskip

For $0 \leq s \leq t \leq T$  let $s = t_0 < t_1 < \dots < t_N = t$ be a temporal partition with uniform stepsize  $t_{n+1} - t_n = \f{t-s}{N} = \Delta t$.   For simplification of notations, we denote
$$\Delta V_{t_n} := V_{t_{n+1}} - V_{t_n},\qquad Y_{n} := Y_{t_{n}}, \qquad Z_{n} := Z_{t_{n}}, \qquad \ora{Y}_n := \ora{Y}_{t_n}, \qquad \ora{Z}_n := \ora{Z}_{t_n}.$$  By Corollary 2.2 in \cite{PP1994}, we have 
$$
\begin{aligned}Y_n(x) &= Y_{t_n}^{t_n, x}, \quad \ora{Y}_{n}(x) &= \ora{Y}_{t_n}^{t_n, x},\quad Y_{n+1}( X_{t_{n+1}}^{t_n, x}) &= Y_{t_{n+1}}^{t_n, x}, \quad \ora{Y}_n(\ola{X}_{t_n}^{t_{n+1}, x}) &= \ora{Y}_{t_n}^{t_{n+1}, x},\\ Z_n(x) &= Z_{t_n}^{t_n, x}, \quad \ora{Z}_{n}(x) &= \ora{Z}_{t_n}^{t_n, x}, \quad Z_{n+1}( X_{t_{n+1}}^{t_n, x}) &= Z_{t_{n+1}}^{t_n, x}, \quad \ora{Z}_n(\ola{X}_{t_n}^{t_{n+1}, x}) &= \ora{Z}_{t_n}^{t_{n+1}, x}.\end{aligned}$$
Denote conditional expectations $$\E[\cdot] := \tilde{\E}[\cdot | \F_T^{V}], \quad \E_x^{n}[\cdot] := \tilde{\E}[\cdot | \mathcal{F}_T^V, X_{t_n} = x  ], \quad  \ola{\E}_x^{n}[\cdot] := \tilde{\E}[\cdot | \mathcal{F}_T^V, \ola{X}_{t_{n}} = x ].$$  It then follows from the definitions of $\E_x^{n}$ and $ \ola{\E}_x^{n}$ that $$\E_x^{n}[Y_n] = Y_n(x), \qquad \ola{\E}_x^{n}[\ora{Y}_n] = \ora{Y}_n(x).$$

Without loss of generality suppose that $\Delta t < s \wedge (T - t) $ and define
\begin{equation*}
Y_N = \f{1}{\Delta t} \int_{t}^{t + \Delta t} Y_r \d r, \hspace{2em}
\ora{Y}_0 = \f{1}{\Delta t} \int_{s - \Delta t}^{s} \ora{Y}_r \d r.
\end{equation*}
For $n = 0, 1, \dots, N-1$, taking the conditional expectations $\E_x^n$ and $\E_x^{\ola{n+1}}$ of temporal discretized approximations of the BDSDEs  in \eqref{Backward:FBDSDEs} and \eqref{Forward:FBDSDEs}, respectively, we have that  (see \cite{Bao})
 \begin{eqnarray}
\E_x^n\left[Y_n\right] &=& \E_x^n\left[Y_{n+1}\right] + \E_x^n\left[h_{n+1} Y_{n+1}\right] \Delta V_{t_n} + \E_x^n\left[\f{\tilde{\rho}_{t_{n+1}}}{\sigma_{t_{n+1}}} Z_{n+1} \right] \Delta V_{t_n}, \\\label{Discrete:Backward}
 \ola{\E}_x^{n+1}[\ora{Y}_{n+1}] &=&  \ola{\E}_x^{n+1}[\ora{Y}_n] +  \ola{\E}_x^{n}\left[ - \ola{b}_n^{\prime} \ora{Y}_n\right] \Delta t \nonumber \\ &&\qquad +  \ola{\E}_x^{n+1}\left[ \ola{h}_n \ora{Y}_n\right] \Delta V_{t_n} -  \ola{\E}_x^{n+1}\left[\f{\tilde{\rho}_{t_{n}}}{\sigma_{t_{n}}} \ora{Z}_{n}\right] \Delta V_{t_n}, \label{Discrete:Forward}
\end{eqnarray}
where  $$h_{n+1} := h( X_{t_{n+1}}), \qquad \ola{b}^{\prime}_n := b_{t_n}^{\prime}(\ola{X}_{t_n}), \qquad \ola{h}_n := h(\ola{X}_{t_n}).$$

By the definition of expectations $\E_x^n$ and $\E_x^{\ola{n+1}}$,
$$\E_x^n\left[ h_{n+1} \right] = \E\left[ h( X_{t_{n+1}}^{t_{n}, x})\right], \quad  \ola{\E}_x^{n+1}\left[ \ola{b}_n^{\prime} \right] = \E\left[ b_n^{\prime}(\ola{X}_{t_n}^{t_{n+1}, x}) \right], \quad \ola{\E}_x^{n+1}\left[ \ola{h}_n \right] = \E\left[ h(\ola{X}_{t_n}^{t_{n+1}, x}) \right].$$
Multiplying  \eqref{Discrete:Backward}  by $\E_x^{\ola{n}}[\ora{Y}_n]$ and \eqref{Discrete:Forward}  by $\E_x^{n+1}[Y_{n+1}]$, then taking integral with respect to $\d x$, we obtain
\begin{equation}\label{Backward:Multiply}
\begin{aligned}
\la \E_x^n\left[Y_n\right] , \ola{\E}_x^{n}[\ora{Y}_n] \ra = & \la \E_x^n\left[Y_{n+1}\right] , \ola{\E}_x^{n}[\ora{Y}_n] \ra  + \la \E_x^n\left[h_{n+1} Y_{n+1}\right] , \ola{\E}_x^{n}[\ora{Y}_n] \ra \Delta V_{t_n} \\
& \quad + \la \E_x^n\left[\f{\tilde{\rho}_{t_{n+1}}}{\sigma_{t_{n+1}}} Z_{n+1} \right] , \ola{\E}_x^{n}[\ora{Y}_n] \ra \Delta V_{t_n}
\end{aligned}
\end{equation}
and
\begin{equation}\label{Forward:Multiply}
\begin{aligned}
& \la \ola{\E}_x^{n+1}[\ora{Y}_{n+1}] , \E_x^{n+1}\left[Y_{n+1}\right] \ra \\
= & \ds \la \ola{\E}_x^{n+1}[\ora{Y}_n] , \E_x^{n+1}\left[Y_{n+1}\right] \ra + \la \ola{\E}_x^{n+1}\left[- \ola{b}_n^{\prime}  \ora{Y}_n \right], \E_x^{n+1}\left[Y_{n+1}\right] \ra  \Delta t \\
& \ds \  + \la \ola{\E}_x^{n+1}\left[\ola{h}_n \ora{Y}_n\right] , \E_x^{n+1}\left[Y_{n+1}\right] \ra \Delta V_{t_n} - \la \ola{\E}_x^{n+1}\left[\f{\tilde{\rho}_{t_{n}}}{\sigma_{t_{n}}} \ora{Z}_{n}\right] , \E_x^{n+1}\left[Y_{n+1}\right] \ra \Delta V_{t_n}.
\end{aligned}
\end{equation}
Subtraction of \eqref{Forward:Multiply} from \eqref{Backward:Multiply} results in
\begin{equation}\label{Subtract}
\begin{aligned}
&  \la \E_x^n\left[Y_n\right] , \ola{\E}_x^{n}[\ora{Y}_n] \ra - \la \ola{\E}_x^{n+1}[\ora{Y}_{n+1}] , \E_x^{n+1}\left[Y_{n+1}\right] \ra \\
 =& \  \underbrace{\la \E_x^{n}\left[Y_{n+1}\right], \ola{\E}_x^{n}[\ora{Y}_n]  - \ola{\E}_x^{n+1}[\ora{Y}_n] \ra + \la \ola{\E}_x^{n+1}[\ora{Y}_n] , \E_x^{n}\left[Y_{n+1}\right]  - \E_x^{n+1}\left[Y_{n+1}\right]  \ra}_{({\bf iv})}\\
 & \ + \underbrace{\la \E_x^n\left[h_{n+1} Y_{n+1}\right] , \ola{\E}_x^{n}[\ora{Y}_n] \ra \Delta V_{t_n} -  \la \ola{\E}_x^{n+1}\left[\ola{h}_n \ora{Y}_n\right] , \E_x^{n+1}\left[Y_{n+1}\right] \ra \Delta V_{t_n}}_{({\bf v})} \\
 & \ + \underbrace{\la \E_x^n\left[\f{\tilde{\rho}_{t_{n+1}}}{\sigma_{t_{n+1}}} Z_{n+1} \right] , \ola{\E}_x^{n}[\ora{Y}_n] \ra \Delta V_{t_n}  + \la \ola{\E}_x^{n+1}\left[\f{\tilde{\rho}_{t_{n}}}{\sigma_{t_{n}}} \ora{Z}_{n}\right] , \E_x^{n+1}\left[Y_{n+1}\right] \ra \Delta V_{t_n}}_{({\bf vi})}\\
  & \ - \la \ola{\E}_x^{n+1}\left[- \ola{b}_n^{\prime}  \ora{Y}_n \right], \E_x^{n+1}\left[Y_{n+1}\right] \ra  \Delta t.
\end{aligned}
\end{equation}
In what follows, we prove that by taking the sum of equation \eqref{Subtract} from $n=0$ to $n = N-1$, the right hand side of the resulting equation converges to $0$ as $\Delta t \to 0$.  To this end we estimate terms ({\bf iv}), ({\bf v}) and ({\bf vi}) one by one.

\smallskip
\noindent ({\bf iv})  By the definitions $\E_x^{\ola{n}}$ and $\E_x^n$, we have
$$\begin{aligned} \ola{\E}_x^{n}[\ora{Y}_n]  - \ola{\E}_x^{n+1}[\ora{Y}_n] &= \E\left[\ora{Y}_{n}(x) - \ora{Y}_{n}( \ola{X}_{t_n}^{t_{n+1, x}} ) \right], \\  \E_x^n\left[Y_{n+1}\right] - \E_x^{n+1}\left[ Y_{n+1} \right] &= \E\left[ Y_{n+1}( X_{t_{n+1}}^{t_n, x}) -  Y_{n+1}(x) \right]. \end{aligned}$$
It follows from It\^o's formula that
\begin{eqnarray}
\ora{Y}_{n}( \ola{X}_{t_n}^{t_{n+1, x}} ) &= & \ora{Y}_{n}(x) +  \int_{t_n}^{t_{n+1}} \big( - b_s(\ola{X}_{s}^{t_{n+1, x}}) \ora{Y}_{n}^{\prime }(\ola{X}_{s}^{t_{n+1, x}})  + \f{(\sigma_s)^2}{2} \ora{Y}_{n}^{ \prime \prime}(\ola{X}_{s}^{t_{n+1, x}}) \big) \d s \nonumber \\
&& \quad + \int_{t_n}^{t_{n+1}} \sigma_s  \ora{Y}_{n}^{\prime}(\ola{X}_{s}^{t_{n+1, x}})  \d \ola{W}_s, \label{estiv-1}\\
 Y_{n+1}(X_{t_{n+1}}^{t_n, x}) &=& \ Y_{n+1}(x) + \int_{t_n}^{t_{n+1}} \big( b_s( X_{s}^{t_n, x} ) Y^{\prime}_{n+1}( X_{s}^{t_n, x} ) + \f{(\sigma_s)^2}{2} Y^{\prime \prime}_{n+1}( X_{s}^{t_n, x} )  \big) \d s \nonumber \\ 
 && \quad + \int_{t_n}^{t_{n+1}} \sigma_s Y^{\prime}_{n+1}( X_{s}^{t_n, x} ) \d W_s. \label{estiv-2}
\end{eqnarray}
Taking conditional expectation $\E$ to Equations \eqref{estiv-1} and \eqref{estiv-2}, we obtain
\begin{eqnarray*}
\ola{\E}_x^{n}[\ora{Y}_n]  - \ola{\E}_x^{n+1}[\ora{Y}_n]  &= & - \E\left[ \int_{t_n}^{t_{n+1}} \big( - b_s(\ola{X}_{s}^{t_{n+1, x}}) \ora{Y}_{n}^{\prime }(\ola{X}_{s}^{t_{n+1, x}})  + \f{(\sigma_s)^2}{2} \ora{Y}_{n}^{ \prime \prime}(\ola{X}_{s}^{t_{n+1, x}}) \big) \d s \right] \\
&= & - \E\left[ - b_n(\ola{X}_{t_n}^{t_{n+1, x}}) \ora{Y}_{n}^{\prime }(\ola{X}_{t_n}^{t_{n+1, x}}) + \f{(\sigma_{t_n})^2}{2} \ora{Y}_{n}^{ \prime \prime}(\ola{X}_{t_{n}}^{t_{n+1, x}})  \right] \cdot \Delta t + \ora{R}_n, \\
\E_x^n\left[Y_{n+1}\right] - \E_x^{n+1}\left[ Y_{n+1} \right] &=&  \E\left[ \int_{t_n}^{t_{n+1}} \big( b_s(  X_{s}^{t_n, x} ) Y^{\prime}_{n+1}( X_{s}^{t_n, x} ) + \f{(\sigma_s)^2}{2} Y^{\prime \prime}_{n+1}( X_{s}^{t_n, x} )  \big) \d s  \right] \nonumber \\
&= & \ \E\left[  b_{n}( X_{t_{n}}^{t_n, x}) Y^{\prime}_{n+1}( X_{t_{n}}^{t_n, x}) + \f{(\sigma_{t_{n}})^2}{2} Y^{\prime \prime}_{n+1}( X_{t_{n}}^{t_n, x} )  \right]  \cdot \Delta t  + R_n,  
\end{eqnarray*}
where
\begin{eqnarray*}
 \ora{R}_n &:= & - \E\left[ \int_{t_n}^{t_{n+1}} \big( - b_s(\ola{X}_{s}^{t_{n+1, x}}) \ora{Y}_{n}^{\prime }(\ola{X}_{s}^{t_{n+1, x}})  + \f{(\sigma_s)^2}{2} \ora{Y}_{n}^{ \prime \prime}(\ola{X}_{s}^{t_{n+1, x}}) \big) \d s \right] \\
 && \ +  \E\left[ - b_n(\ola{X}_{t_n}^{t_{n+1, x}}) \ora{Y}_{n}^{\prime }(\ola{X}_{t_n}^{t_{n+1, x}}) + \f{(\sigma_{t_n})^2}{2} \ora{Y}_{n}^{ \prime \prime}(\ola{X}_{t_{n}}^{t_{n+1, x}}) \right]  \cdot \Delta t, \\
R_n &:=&  \E\left[ \int_{t_n}^{t_{n+1}} \big( b_s( X_{s}^{t_n, x} ) Y^{\prime}_{n+1}( X_{s}^{t_n, x} ) + \f{(\sigma_s)^2}{2} Y^{\prime \prime}_{n+1}( X_{s}^{t_n, x} )  \big) \d s  \right] \\
&& \ - \E\big[  b_{n}( X_{t_{n}}^{t_n, x}) Y^{\prime}_{n+1}( X_{t_{n}}^{t_n, x}) + \f{(\sigma_{t_{n}})^2}{2} Y^{\prime \prime}_{n+1}( X_{t_{n}}^{t_n, x} )  \big]  \cdot \Delta t .
\end{eqnarray*}
As a consequence 
\begin{equation}\label{eq:Part1:a}
\begin{aligned}
&\la \E_x^{n}\left[Y_{n+1}\right], \ola{\E}_x^{n}[\ora{Y}_n]  - \ola{\E}_x^{n+1}[\ora{Y}_n] \ra \\
= & \  - \int_{\mathbb{R}} \E\left[Y_{n+1}( X_{t_{n+1}}^{t_n, x})\right]   \Big( \E\big[  - b_n(\ola{X}_{t_n}^{t_{n+1, x}}) \ora{Y}_{n}^{\prime }(\ola{X}_{t_n}^{t_{n+1, x}}) + \f{(\sigma_{t_n})^2}{2} \ora{Y}_{n}^{ \prime \prime}(\ola{X}_{t_{n}}^{t_{n+1, x}}) \big]  \cdot \Delta t - \ora{R}_n \Big)  \d x.  
\end{aligned}
\end{equation}
Similarly 
\begin{equation}\label{eq:Part1:b}
\begin{aligned}
& \la \ola{\E}_x^{n+1}[\ora{Y}_n] , \E_x^{n}\left[Y_{n+1}\right]  - \E_x^{n+1}\left[Y_{n+1}\right]  \ra\\
= & \ \int_{\mathbb{R}} \E\left[\ora{Y}_{n}(\ola{X}_{t_n}^{t_{n+1, x}})\right] \Big( \E\big[  b_{n}( X_{t_{n}}^{t_n, x}) Y^{\prime}_{n+1}( X_{t_{n}}^{t_n, x}) + \f{(\sigma_{t_{n}})^2}{2} Y^{\prime \prime}_{n+1}( X_{t_{n}}^{t_n, x} )  \big]  \cdot \Delta t + R_n \Big)  \d x.
\end{aligned}
\end{equation}
Adding \eqref{eq:Part1:a} to \eqref{eq:Part1:b} we have that 
\begin{equation}\label{eq:Sum:main:total}
\begin{aligned}
({\bf iv}) &=  \ \Big( \ \ \underbrace{- \int_{\mathbb{R}} \E\left[Y_{n+1}( X_{t_{n+1}}^{t_n, x})\right]   \E\big[  - b_n(\ola{X}_{t_n}^{t_{n+1, x}}) \ora{Y}_{n}^{\prime }(\ola{X}_{t_n}^{t_{n+1, x}})  \big] \d x}_{({\bf iv-1})} \\
& \qquad + \underbrace{\int_{\mathbb{R}} \E\left[\ora{Y}_{n}(\ola{X}_{t_n}^{t_{n+1, x}})\right]  \E\big[  b_{n}( X_{t_{n}}^{t_n, x}) Y^{\prime}_{n+1}( X_{t_{n}}^{t_n, x}) \big] \d x }_{({\bf iv-2})} \ \ \Big)  \cdot \Delta t\\
& \qquad +\Big( \ \  \underbrace{ - \int_{\mathbb{R}} \E\left[Y_{n+1}( X_{t_{n+1}}^{t_n, x})\right]  \E\big[ \f{(\sigma_{t_n})^2}{2} \ora{Y}_{n}^{ \prime \prime}(\ola{X}_{t_{n}}^{t_{n+1, x}})  \big] \d x}_{({\bf iv-3})} \\
& \qquad +  \underbrace{\int_{\mathbb{R}} \E\left[\ora{Y}_{n}(\ola{X}_{t_n}^{t_{n+1, x}})\right]  \E\big[ \f{(\sigma_{t_{n}})^2}{2} Y^{\prime \prime}_{n+1}( X_{t_{n}}^{t_n, x} )  \big] \d x }_{({\bf iv-4})} \ \ \Big) \cdot \Delta t  + R_n^x,
\end{aligned}
\end{equation}
where
$$R_n^x = \int_{\mathbb{R}} \E[Y_{n+1}( X_{t_{n+1}}^{t_n, x})]  \ora{R}_n dx + \int_{\mathbb{R}} \E[\ora{Y}_{n}(\ola{X}_{t_n}^{t_{n+1, x}})]  R_n \d x. $$
Again by using the It\^o formula we obtain
\begin{eqnarray*}
\ora{Y}^{\prime}_{n}( \ola{X}_{t_n}^{t_{n+1, x}} ) &= & \ora{Y}^{\prime}_{n}(x) +  \int_{t_n}^{t_{n+1}} \big( - b_s(\ola{X}_{s}^{t_{n+1, x}}) \ora{Y}_{n}^{ \prime \prime }(\ola{X}_{s}^{t_{n+1, x}})  + \f{(\sigma_s)^2}{2} \ora{Y}_{n}^{(3)}(\ola{X}_{s}^{t_{n+1, x}}) \big) \d s \\
&& \qquad + \int_{t_n}^{t_{n+1}} \sigma_s  \ora{Y}_{n}^{\prime \prime}(\ola{X}_{s}^{t_{n+1, x}})  \d\ola{W}_s, \\
\ora{Y}_n^{\prime \prime}(\ola{X}_{t_n}^{t_{n+1},x}) &= & \ora{Y}^{\prime \prime}_{n}(x) +  \int_{t_n}^{t_{n+1}} \big( - b_s(\ola{X}_{s}^{t_{n+1, x}}) \ora{Y}_{n}^{ (3) }(\ola{X}_{s}^{t_{n+1, x}})  + \f{(\sigma_s)^2}{2} \ora{Y}_{n}^{(4)}(\ola{X}_{s}^{t_{n+1, x}}) \big) \d s \\
&& \qquad + \int_{t_n}^{t_{n+1}} \sigma_s  \ora{Y}_{n}^{(3)}(\ola{X}_{s}^{t_{n+1, x}})  \d \ola{W}_s, \\
b_n(\ola{X}_{t_n}^{t_{n+1, x}}) &= & b_n(x) + \int_{t_n}^{t_{n+1}} - b_s(\ola{X}_{s}^{t_{n+1, x}}) b_n^{\prime}(\ola{X}_{s}^{t_{n+1, x}})  + \f{(\sigma_s)^2}{2} b_{n}^{\prime \prime}(\ola{X}_{s}^{t_{n+1, x}}) \big) \d s \\
&& \qquad + \int_{t_n}^{t_{n+1}} \sigma_s  b_{n}^{ \prime}(\ola{X}_{s}^{t_{n+1, x}})  \d \ola{W}_s.
\end{eqnarray*}
Hence, the term $\E\big[  b_n(\ola{X}_{t_n}^{t_{n+1, x}}) \ora{Y}_{n}^{\prime }(\ola{X}_{t_n}^{t_{n+1, x}})  \big]$ on the right hand side of \eqref{eq:Sum:main:total} can be written as
$
 \E\big[  b_n(\ola{X}_{t_n}^{t_{n+1, x}}) \ora{Y}_{n}^{\prime }(\ola{X}_{t_n}^{t_{n+1, x}})  \big] =  b_n(x) \ora{Y}_n^{\prime}(x) + P_n(x)$ with
$$ 
\begin{aligned}
P_n(x) & = \E\Big[ b_n(x) \cdot \int_{t_n}^{t_{n+1}} \big( - b_s(\ola{X}_{s}^{t_{n+1, x}}) \ora{Y}_{n}^{ \prime \prime }(\ola{X}_{s}^{t_{n+1, x}})  + \f{(\sigma_s)^2}{2} \ora{Y}_{n}^{(3)}(\ola{X}_{s}^{t_{n+1, x}}) \big) \d s \\
& \quad + \ora{Y}^{\prime}_{n}(x) \cdot \int_{t_n}^{t_{n+1}} \big( - b_s(\ola{X}_{s}^{t_{n+1, x}}) b_n^{\prime}(\ola{X}_{s}^{t_{n+1, x}})  + \f{(\sigma_s)^2}{2} b_{n}^{\prime \prime}(\ola{X}_{s}^{t_{n+1, x}}) \big) \d s   \\
& \quad + \int_{t_n}^{t_{n+1}} \big( - b_s (\ola{X}_{s}^{t_{n+1, x}}) b_n^{\prime}(\ola{X}_{s}^{t_{n+1, x}})  + \f{(\sigma_s)^2}{2} b_{n}^{\prime \prime}(\ola{X}_{s}^{t_{n+1, x}}) \big) \d s \\
& \qquad  \cdot \int_{t_n}^{t_{n+1}} \big( - b_s(\ola{X}_{s}^{t_{n+1, x}}) \ora{Y}_{n}^{ \prime \prime }(\ola{X}_{s}^{t_{n+1, x}})  + \f{(\sigma_s)^2}{2} \ora{Y}_{n}^{(3)}(\ola{X}_{s}^{t_{n+1, x}}) \big) \d s  \\
& \quad  + \int_{t_n}^{t_{n+1}} (\sigma_s)^2 b_n^{\prime}(\ola{X}_{s}^{t_{n+1, x}}) \ora{Y}_n^{\prime \prime}(\ola{X}_{s}^{t_{n+1, x}}) \d s  \Big]. 
\end{aligned}
$$
As a  result the terms on the right hand side of \eqref{eq:Sum:main:total} can be rewritten as 
\begin{eqnarray}
({\bf iv-1}) &=& \int_{\mathbb{R}} \big(  Y_{n+1}(x) \cdot b_n(x) \ora{Y}_n^{\prime}(x) \big) \d x + H_n^{1}, \label{eq:Sum:main:1} \\
({\bf iv-2}) &=& \int_{\mathbb{R}} \big(  \ora{Y}_n(x) b_{n}(x) Y_{n+1}^{\prime}(x) \big) \d x + H_n^2, \label{eq:Sum:main:2} \\
({\bf iv-3}) &=& - \int_{\mathbb{R}}   \f{(\sigma_{t_n})^2}{2}  Y_{n+1}(x)\ora{Y}^{\prime \prime}_{n}(x) \d x + H_n^3, \label{eq:Sum:main:3} \\
({\bf iv-4}) &=&  \int_{\mathbb{R}} \f{(\sigma_{t_{n}})^2}{2} \ora{Y}_{n}(x)  Y_{n+1}^{\prime \prime}(x)  \d x + H_n^4,\label{eq:Sum:main:4}
\end{eqnarray}
where
$$
\begin{aligned}
H_n^1 =&  \int_{\mathbb{R}} \big\{ \Big( Y_{n+1}(x) + \E[ \int_{t_n}^{t_{n+1}} \big( b_s(  X_{s}^{t_n, x} ) Y^{\prime}_{n+1}( X_{s}^{t_n, x} ) + \f{(\sigma_s)^2}{2} Y^{\prime \prime}_{n+1}( X_{s}^{t_n, x} )  \big) ds  \big)] \Big) \cdot  P_n(x)  \\
& \quad  +  \E[ \int_{t_n}^{t_{n+1}} \big( b_s(  X_{s}^{t_n, x} ) Y^{\prime}_{n+1}( X_{s}^{t_n, x} ) + \f{(\sigma_s)^2}{2} Y^{\prime \prime}_{n+1}( X_{s}^{t_n, x} )  \big) ds  \big)] \cdot  b_n(x) \ora{Y}_n^{\prime}(x) \big\} \d x, \\
H_n^2 =&  \int_{\mathbb{R}} \big\{ \E[ \int_{t_n}^{t_{n+1}} \big( - b_s(\ola{X}_{s}^{t_{n+1, x}}) \ora{Y}_{n}^{\prime }(\ola{X}_{s}^{t_{n+1, x}})  + \f{(\sigma_s)^2}{2} \ora{Y}_{n}^{ \prime \prime}(\ola{X}_{s}^{t_{n+1, x}}) \big) ds ] \cdot   b_{n}(x) Y_{n+1}^{\prime}(x) \big\} \d x  , \\
H_n^3 =&   - \int_{\mathbb{R}} \big\{ Y_{n+1}(x) 
\cdot \f{(\sigma_{t_n})^2}{2} \E[  \int_{t_n}^{t_{n+1}} \big( - b_s(\ola{X}_{s}^{t_{n+1, x}}) \ora{Y}_{n}^{ (3) }(\ola{X}_{s}^{t_{n+1, x}})  + \f{(\sigma_s)^2}{2} \ora{Y}_{n}^{(4)}(\ola{X}_{s}^{t_{n+1, x}}) \big) \d s  ] \\
& \quad + \E[ \int_{t_n}^{t_{n+1}} \big( b_s( X_{s}^{t_n, x} ) Y^{\prime}_{n+1}( X_{s}^{t_n, x} ) + \f{(\sigma_s)^2}{2} Y^{\prime \prime}_{n+1}( X_{s}^{t_n, x} )  \big) \d s  \big)]   \cdot \f{(\sigma_{t_n})^2}{2} \ora{Y}^{\prime \prime}_{n}(x) \\
& \quad + \E[ \int_{t_n}^{t_{n+1}} \big( b_s(  X_{s}^{t_n, x} ) Y^{\prime}_{n+1}( X_{s}^{t_n, x} ) + \f{(\sigma_s)^2}{2} Y^{\prime \prime}_{n+1}( X_{s}^{t_n, x} )  \big) \d s  \big)]  \\
& \qquad \cdot  \f{(\sigma_{t_n})^2}{2} \E[  \int_{t_n}^{t_{n+1}} \big( - b_s(\ola{X}_{s}^{t_{n+1, x}}) \ora{Y}_{n}^{ (3) }(\ola{X}_{s}^{t_{n+1, x}})  + \f{(\sigma_s)^2}{2} \ora{Y}_{n}^{(4)}(\ola{X}_{s}^{t_{n+1, x}}) \big) \d s  ]  \big\}\d x, \\
H_n^4 =& \int_{\mathbb{R}} \big\{ \E[ \int_{t_n}^{t_{n+1}} \big( - b_s(\ola{X}_{s}^{t_{n+1, x}}) \ora{Y}_{n}^{\prime }(\ola{X}_{s}^{t_{n+1, x}})  + \f{(\sigma_s)^2}{2} \ora{Y}_{n}^{ \prime \prime}(\ola{X}_{s}^{t_{n+1, x}}) \big)] \d s 
\cdot \f{(\sigma_{t_{n}})^2}{2}  Y_{n+1}^{\prime \prime}(x) \big\} \d x. 
\end{aligned}$$
\smallskip

Integrating by parts, we obtain 
\begin{eqnarray}
 \int_{\mathbb{R}} \big(  Y_{n+1}(x) \cdot b_n(x) \ora{Y}_n^{\prime}(x) \big) \d x &=& - \int_{\mathbb{R}} Y_{n+1}^{\prime}(x) b_n(x) \ora{Y}_n(x)  \d x \nonumber \\ && - \int_{\mathbb{R}} Y_{n+1}(x) b_n^{\prime}(x) \ora{Y}_n(x)  \d x,  \label{Integration-by-parts:1} \\
 - \int_{\mathbb{R}}   \f{(\sigma_{t_n})^2}{2}  Y_{n+1}(x)\ora{Y}^{\prime \prime}_{n}(x) \d x &=& \int_{\mathbb{R}}   \f{(\sigma_{t_n})^2}{2}  Y_{n+1}^{\prime}(x)\ora{Y}^{\prime}_{n}(x) \d x, \\ \label{Integration-by-parts:2}
\int_{\mathbb{R}} \f{(\sigma_{t_{n}})^2}{2} \ora{Y}_{n}(x)  Y_{n+1}^{\prime \prime}(x)  \d x  &=& - \int_{\mathbb{R}} \f{(\sigma_{t_{n}})^2}{2} \ora{Y}_{n}^{\prime}(x)  Y_{n+1}^{\prime }(x)  \d x. \label{Integration-by-parts:3}
\end{eqnarray}
Adding \eqref{eq:Sum:main:1} to \eqref{eq:Sum:main:2} and applying \eqref{Integration-by-parts:1}, the sum of the first two terms on the right hand side of \eqref{eq:Sum:main:total} becomes
\begin{equation}\label{eq:Sum:main:1+2}
\begin{aligned}
&  - \int_{\mathbb{R}} \E[Y_{n+1}( X_{t_{n+1}}^{t_n, x})]   \E\big[  - b_n(\ola{X}_{t_n}^{t_{n+1, x}}) \ora{Y}_{n}^{\prime }(\ola{X}_{t_n}^{t_{n+1, x}})  \big] dx \\
& \ \ + \int_{\mathbb{R}} \E[\ora{Y}_{n}(\ola{X}_{t_n}^{t_{n+1, x}})]  \E\big[  b_{n}( X_{t_{n}}^{t_n, x}) Y^{\prime}_{n+1}( X_{t_{n}}^{t_n, x}) \big] dx \\
= &  - \int_{\mathbb{R}} Y_{n+1}(x) b_n^{\prime}(x) \ora{Y}_n(x)  dx   +H_n^1 + H_n^2.  
\end{aligned}
\end{equation}
Similarly, adding \eqref{eq:Sum:main:3} to \eqref{eq:Sum:main:4} and applying \eqref{Integration-by-parts:2}, \eqref{Integration-by-parts:3} yields 
\begin{equation}\label{eq:Sum:main:3+4}
\begin{aligned}
& - \int_{\mathbb{R}} \E[Y_{n+1}(\tilde{X}_{t_{n+1}}^{t_n, x})]  \E\big[ \f{(\sigma_{t_n})^2}{2} \ora{Y}_{n}^{ \prime \prime}(\ola{X}_{t_{n}}^{t_{n+1, x}})  \big] dx + \int_{\mathbb{R}} \E[\ora{Y}_{n}(\ola{X}_{t_n}^{t_{n+1, x}})]  \E\big[ \f{(\sigma_{t_{n}})^2}{2} Y^{\prime \prime}_{n+1}(\tilde{X}_{t_{n}}^{t_n, x} )  \big] dx \\
& = H_n^3 + H_n^4,
\end{aligned}
\end{equation}
which is the sum of the last two terms on the right hand side of \eqref{eq:Sum:main:total}.

For first two terms on the right hand side of equation \eqref{Subtract},  we insert \eqref{eq:Sum:main:1+2} and \eqref{eq:Sum:main:3+4} into \eqref{eq:Sum:main:total} to obtain the following equation
\begin{equation}\label{eq:Sum:main}
\begin{aligned}
&\la \E_x^{n}[Y_{n+1}], \ola{\E}_x^{n}[\ora{Y}_n]  - \ola{\E}_x^{n+1}[\ora{Y}_n] \ra  + \la \ola{\E}_x^{n+1}[\ora{Y}_n] , \E_x^{n}[Y_{n+1}]  - \E_x^{n+1}[Y_{n+1}]  \ra\\
= & \  \big( - \int_{\mathbb{R}} Y_{n+1}(x) b_n^{\prime}(x) \ora{Y}_n(x)  dx \big) \Delta t   + ( H_n^1 + H_n^2 + H_n^3 + H_n^4 ) \Delta t + R_n^x.  
\end{aligned}
\end{equation}

Next, we consider the term 
$$
\la \E_x^n[h_{n+1} Y_{n+1}] , \ola{\E}_x^{n}[\ora{Y}_n] \ra \Delta V_{t_n} -  \la \ola{\E}_x^{n+1}[\ola{h}_n \ora{Y}_n] , \E_x^{n+1}[Y_{n+1}] \ra \Delta V_{t_n}
$$ in \eqref{Subtract}. From the definition of $\E_x^n$ and $\ola{\E}_x^{n+1}$, one has
\begin{equation}\label{eq:sum-part2-1}
\begin{aligned}
&  \la \E_x^n[h_{n+1} Y_{n+1}] , \ola{\E}_x^{n}[\ora{Y}_n] \ra \Delta V_{t_n} -  \la \ola{\E}_x^{n+1}[\ola{h}_n \ora{Y}_n] , \E_x^{n+1}[Y_{n+1}] \ra \Delta V_{t_n} \\
= & \ \int_{\mathbb{R}} \E[ h( X_{t_{n+1}}^{t_n, x})  Y_{n+1}( X_{t_{n+1}}^{t_n, x}) \ora{Y}_n(x) ] - \E[ h(\ola{X}_{t_n}^{t_{n+1}, x}) \ora{Y}_{n}(\ola{X}_{t_n}^{t_{n+1}, x})  Y_{n+1}(x) ] dx  \Delta V_{t_n}.
\end{aligned}
\end{equation}
We apply It\^o formula to $h$ on time interval $[t_n, t_{n+1}]$ to get
$$
\begin{aligned}
h( X_{t_{n+1}}^{t_n, x}) = &h(x) + \int_{t_n}^{t_{n+1}} b_s ( X_{s}^{t_n, x}) h^{\prime}( X_{s}^{t_n, x})  + \f{(\sigma_s)^2}{2} h^{\prime \prime}( X_{s}^{t_n, x}) \big) ds \\
& \ + \int_{t_n}^{t_{n+1}} \sigma_s  h^{ \prime}( X_{s}^{t_n, x})  d W_s,
\end{aligned}
$$
and
$$
\begin{aligned}
h(\ola{X}_{t_n}^{t_{n+1, x}}) = & h(x) + \int_{t_n}^{t_{n+1}} - b_s (\ola{X}_{s}^{t_{n+1, x}}) h^{\prime}(\ola{X}_{s}^{t_{n+1, x}})  + \f{(\sigma_s)^2}{2} h^{\prime \prime}(\ola{X}_{s}^{t_{n+1, x}}) \big) ds \\
& \ + \int_{t_n}^{t_{n+1}} \sigma_s  h^{ \prime}(\ola{X}_{s}^{t_{n+1, x}})  d\ola{W}_s.
\end{aligned}
$$
Thus 
$$
\begin{aligned}
& \E[ h( X_{t_{n+1}}^{t_n, x})  Y_{n+1}( X_{t_{n+1}}^{t_n, x}) \ora{Y}_n(x) ] \\
= & \ h(x) Y_{n+1}(x) \ora{Y}_n(x) + \ora{Y}_n(x) \cdot \E\big[ \big( h( X_{t_{n+1}}^{t_n, x}) - h(x) \big)  Y_{n+1}(x) \\
& \quad + h(x) \big( Y_{n+1}( X_{t_{n+1}}^{t_n, x}) - Y_{n+1}(x) \big) + \big( h( X_{t_{n+1}}^{t_n, x}) - h(x) \big) \big( Y_{n+1}( X_{t_{n+1}}^{t_n, x}) - Y_{n+1}(x) \big)  \big]
\end{aligned}
$$
and
$$
\begin{aligned}
& \E[ h(\ola{X}_{t_n}^{t_{n+1}, x}) \ora{Y}_{n}(\ola{X}_{t_n}^{t_{n+1}, x})  Y_{n+1}(x) ]  \\
= & \ h(x) \ora{Y}_n(x)  Y_{n+1}(x) + Y_{n+1}(x) \cdot \E\big[ \big( h(\ola{X}_{t_{n}}^{t_{n+1}, x}) - h(x) \big)  \ora{Y}_{n}(x) \\
& \quad + h(x) \big( \ora{Y}_{n}(\ola{X}_{t_{n}}^{t_{n+1}, x}) - \ora{Y}_{n}(x) \big) + \big( h(\ola{X}_{t_{n}}^{t_{n+1}, x}) - h(x) \big)    \big( \ora{Y}_{n}(\ola{X}_{t_{n}}^{t_{n+1}, x}) - \ora{Y}_{n}(x) \big) \big].
\end{aligned}
$$
With the above equations, \eqref{eq:sum-part2-1} becomes
\begin{equation}\label{eq:sum-part2-2}
\begin{aligned}
  \la \E_x^n[h_{n+1} Y_{n+1}] , \ola{\E}_x^{n}[\ora{Y}_n] \ra \Delta V_{t_n} -  \la \ola{\E}_x^{n+1}[\ola{h}_n \ora{Y}_n] , \E_x^{n+1}[Y_{n+1}] \ra \Delta V_{t_n} =  G_n^1 \Delta V_{t_n},
\end{aligned}
\end{equation}
where
$$
\begin{aligned}
G_n^1 = & \ \int_{\mathbb{R}} \big\{ \ora{Y}_n(x) \cdot \E\big[ \big( h( X_{t_{n+1}}^{t_n, x}) - h(x) \big)  Y_{n+1}(x) \\
& \qquad + h(x) \big( Y_{n+1}( X_{t_{n+1}}^{t_n, x}) - Y_{n+1}(x) \big) + \big( h( X_{t_{n+1}}^{t_n, x}) - h(x) \big) \big( Y_{n+1}( X_{t_{n+1}}^{t_n, x}) - Y_{n+1}(x) \big)  \big] \\
& \ - Y_{n+1}(x) \cdot \E\big[ \big( h(\ola{X}_{t_{n}}^{t_{n+1}, x}) - h(x) \big)  \ora{Y}_{n}(x) \\
& \qquad + h(x) \big( \ora{Y}_{n}(\ola{X}_{t_{n}}^{t_{n+1}, x}) - \ora{Y}_{n}(x) \big) + \big( h(\ola{X}_{t_{n}}^{t_{n+1}, x}) - h(x) \big)    \big( \ora{Y}_{n}(\ola{X}_{t_{n}}^{t_{n+1}, x}) - \ora{Y}_{n}(x) \big) \big] \big\} dx .
\end{aligned}
$$

Finally, we consider the term 
$$\la \E_x^n[\f{\tilde{\rho}_{t_{n+1}}}{\sigma_{t_{n+1}}} Z_{n+1} ] , \ola{\E}_x^{n}[\ora{Y}_n] \ra \Delta V_{t_n}  + \la \ola{\E}_x^{n+1}[\f{\tilde{\rho}_{t_{n}}}{\sigma_{t_{n}}} \ora{Z}_{n}] , \E_x^{n+1}[Y_{n+1}] \ra \Delta V_{t_n}$$ 
on the right hand side of equation \eqref{Subtract}. From the relation between $Z_t$ and $\f{\p Y_t}{\p x}$ given in \eqref{Z=dY}, we know that 
$$Z_{n+1}(X_{t_{n+1}}^{t_n, x}) = \f{\p Y_{n+1}(X_{t_{n+1}}^{t_n, x})}{\p x} (\nabla X_{t_{n+1}}^{t_n, x})^{-1} \sigma_{t_{n+1}}$$
and 
$$\ora{Z}_{n}(\ola{X}_{t_n}^{t_{n+1}, x}) = \f{\p \ora{Y}_{n}(\ola{X}_{t_n}^{t_{n+1}, x}) }{\p x}(\nabla \ola{X}_{t_{n}}^{t_{n+1}, x})^{-1} \sigma_{t_{n}}.$$
Therefore we have
\begin{equation}\label{Z-Y:1}
\begin{aligned}
\la \E_x^n[\f{\tilde{\rho}_{t_{n+1}}}{\sigma_{t_{n+1}}} Z_{n+1} ] , \ola{\E}_x^{n}[\ora{Y}_n] \ra
= & \ \E \Big[ \int_{\mathbb{R}} \f{\tilde{\rho}_{t_{n+1}}}{\sigma_{t_{n+1}}} Z_{n+1} (X_{t_{n+1}}^{t_n, x}) \cdot \ora{Y}_{n} (x) dx \Big]\\
= & \ \E \Big[ \int_{\mathbb{R}} \tilde{\rho}_{t_{n+1}} \f{\p Y_{n+1}(X_{t_{n+1}}^{t_n, x}) }{\p x} \cdot \ora{Y}_{n} (x) \cdot (\nabla X_{t_{n+1}}^{t_n, x})^{-1}  dx \Big]
\end{aligned}
\end{equation}
and
\begin{equation}\label{Z-Y:2}
\begin{aligned}
 \la \ola{\E}_x^{n+1}[\f{\tilde{\rho}_{t_{n}}}{\sigma_{t_{n}}} \ora{Z}_{n}] , \E_x^{n+1}[Y_{n+1}] \ra 
= & \  \E \Big[ \int_{\mathbb{R}} \f{\tilde{\rho}_{t_{n}}}{\sigma_{t_{n}}} \ora{Z}_{n} (\ola{X}_{t_n}^{t_{n+1}, x}) \cdot Y_{n+1} (x) dx \Big] \\
= & \ \E \Big[ \int_{\mathbb{R}} \tilde{\rho}_{t_{n}} \f{\p \ora{Y}_{n}(\ola{X}_{t_n}^{t_{n+1}, x})}{\p x}  \cdot Y_{n+1} (x) \cdot (\nabla \ola{X}_{t_{n}}^{t_{n+1}, x})^{-1} dx \Big],
\end{aligned}
\end{equation}
Adding \eqref{Z-Y:1} and \eqref{Z-Y:2} together, we obtain
\begin{equation}\label{Z-Y}
\begin{aligned}
 & \la \E_x^n[\f{\tilde{\rho}_{t_{n+1}}}{\sigma_{t_{n+1}}} Z_{n+1} ] , \ola{\E}_x^{n}[\ora{Y}_n] \ra + \la \ola{\E}_x^{n+1}[\f{\tilde{\rho}_{t_{n}}}{\sigma_{t_{n}}} \ora{Z}_{n}] , \E_x^{n+1}[Y_{n+1}] \ra  \\
 = & \ \E \Big[ \int_{\mathbb{R}} \tilde{\rho}_{t_{n+1}} \f{\p Y_{n+1} }{\p x}(x) \cdot \ora{Y}_{n} (x)  dx  + \int_{\mathbb{R}} \tilde{\rho}_{t_{n+1}} \f{\p \ora{Y}_{n} }{\p x}(x)  \cdot Y_{n+1} (x) dx \Big] + G_n^2,
\end{aligned}
\end{equation}
where 
$$
\begin{aligned}
G_n^2 &= \E \Big[ \int_{\mathbb{R}} \tilde{\rho}_{t_{n+1}} \f{\p Y_{n+1} }{\p x}(X_{t_{n+1}}^{t_n, x}) \cdot \ora{Y}_{n} (x) \cdot (\nabla X_{t_{n+1}}^{t_n, x})^{-1}  dx  - \int_{\mathbb{R}} \tilde{\rho}_{t_{n+1}} \f{\p Y_{n+1} }{\p x}(x) \cdot \ora{Y}_{t_{n}} (x)  dx \Big] \\
& \ + \E \Big[ \int_{\mathbb{R}} \tilde{\rho}_{t_{n}} \f{\p \ora{Y}_{n} }{\p x} (\ola{X}_{t_n}^{t_{n+1}, x})) \cdot Y_{n+1} (x) \cdot (\nabla \ola{X}_{t_{n}}^{t_{n+1}, x})^{-1} dx - \int_{\mathbb{R}} \tilde{\rho}_{t_{n+1}} \f{\p \ora{Y}_{n} }{\p x}(x)  \cdot Y_{n+1} (x) dx \Big].
\end{aligned}
$$
Integrating  by parts gives, 
$$ \int_{\mathbb{R}} \tilde{\rho}_{t_{n+1}} \f{\p Y_{n+1} }{\p x}(x) \cdot \ora{Y}_{n} (x)  dx =  - \int_{\mathbb{R}} \tilde{\rho}_{t_{n+1}} \f{\p \ora{Y}_{n} }{\p x}(x)  \cdot Y_{n+1} (x) dx.$$
Therefore 
\begin{equation}\label{eq:Z-Y}
\begin{aligned}
G_n^2=& \la \E_x^n[\f{\tilde{\rho}_{t_{n+1}}}{\sigma_{t_{n+1}}} Z_{n+1} ] , \E_x^{\ola{n}}[\ora{Y}_n] \ra \Delta V_{t_n}  + \la \E_x^{\ola{n+1}}[\f{\tilde{\rho}_{t_{n}}}{\sigma_{t_{n}}} \ora{Z}_{n}] , \E_x^{n+1}[Y_{n+1}] \ra. 
\end{aligned}
\end{equation}

From \eqref{eq:Sum:main}, \eqref{eq:sum-part2-2} and \eqref{eq:Z-Y}, equation \eqref{Subtract} becomes
\begin{equation}\label{eq:sum:local}
\begin{aligned}
&  \la \E_x^n[Y_n] , \ola{\E}_x^{n}[\ora{Y}_n] \ra - \la \E_x^{n+1}[Y_{n+1}] , \ola{\E}_x^{n+1}[\ora{Y}_{n+1}] \ra \\
 = &  \big( - \int_{\mathbb{R}} Y_{n+1}(x) b_n^{\prime}(x) \ora{Y}_n(x)  dx \big) \Delta t   + ( H_n^1 + H_n^2 + H_n^3 + H_n^4 ) \Delta t + R_n^x \\
 & \ + \int_{\mathbb{R}} \E[ b_n^{\prime}(\ola{X}_{t_{n}}^{t_{n+1}, x}) \ora{Y}_n(\ola{X}_{t_{n}}^{t_{n+1}, x})] Y_{n+1}(x)  dx \Delta t  + G_n  \Delta V_{t_n} \\
 = & H_n \Delta t + R_n^x + G_n \Delta V_{t_n} + F_n \Delta t, 
\end{aligned}
\end{equation}
where 
$$H_n = H_n^1 + H_n^2 + H_n^3 + H_n^4, \quad G_n = G_n^1 + G_n^2$$
and
$$F_n = \int_{\mathbb{R}}Y_{n+1}(x)  \Big(  \E\big[ b_n^{\prime}(\ola{X}_{t_{n}}^{t_{n+1}, x}) \ora{Y}_n(\ola{X}_{t_{n}}^{t_{n+1}, x})\big] -  b_n^{\prime}(x) \ora{Y}_n(x) \Big) dx .$$
Next,  we sum \eqref{eq:sum:local} from $n = 0$ to $n = N-1$ to get
\begin{equation}\label{eq:sum:global}
\begin{aligned}
& \la \E_x^0[Y_0] , \ola{\E}_x^{0}[\ora{Y}_0] \ra - \la \ola{\E}_x^{N}[\ora{Y}_{N}] , \E_x^{N}[Y_{N}] \ra \\
= & \sum_{n=0}^{N-1} ( H_n \Delta t + R_n^x + G_n \Delta V_{t_n} + F_n \Delta t ) .
\end{aligned}
\end{equation}
From definitions of $H_n$ , $R_n^x$, $G_n$ and $F_n$, it's easy to verify that
$E[(H_n)^2] \leq C (\Delta t)^2$,
$E[(R_n^x)^2] \leq C (\Delta t)^{4}$,$
E[(G_n)^2] \leq C (\Delta t)^{2}$
and $E[(F_n)^2] \leq C(\Delta t)$.
Therefore,
$$ \lim_{\Delta t \rightarrow 0} \sum_{n=0}^{N-1} ( H_n \Delta t + R_n^x + G_n \Delta V_{t_n} + F_n \Delta t ) = 0, \ a.s. .$$
Also, since $\lim_{\Delta t \rightarrow  0}Y_0 = Y_s$ and $\lim_{\Delta t \rightarrow  0}Y_N = Y_t$,
we have
\begin{equation*}
\la Y_s, \ora{Y}_s \ra = \la Y_t, \ora{Y}_t \ra
\end{equation*}
as required.
$\Box$

\bigskip
Now are ready to state the main result in this paper. It is a direct consequence of  Theorems \ref{thm:Feynman-Kac} and \ref{thm:Adjoint-Operator}.   
\begin{thm}\label{Thm:Final}
Assume that the ssumptions in Theorem \ref{thm:Feynman-Kac} and Theorem \ref{thm:Adjoint-Operator} hold. Then 
$$\la \ora{Y}_T, \phi \ra = \tilde{E}\left[\phi(U_T)Q_T \big| \mathcal{F}_T^V\right], \hspace{2em} \forall \phi \in \L^{\infty}(\mathbb{R}^d) . $$ 
\end{thm}
\textbf{Proof.}
Applying Theorem \ref{thm:Adjoint-Operator}, one has
$$ \la\ora{Y}_T, Y_T \ra = \la \ora{Y}_0, Y_0\ra. $$
Since $Y_T = \phi$ as given in \eqref{Backward:FBDSDEs}, $\ora{Y}_0 = p_0$ as given in \eqref{Forward:FBDSDEs} and $Y_0 = \tilde{E}_x\left[\phi(S_T) Q_T \big| \mathcal{F}_T^V \right]$ as proved in Theorem \ref{thm:Feynman-Kac}, we have
$$\la \ora{Y}_T, \phi \ra = \int_{\mathbb{R}} p_0(x) \tilde{E}_x[ \phi(U_T) Q_T \big| \mathcal{F}_T^V ] dx.  $$
Let $\varphi $ be any bounded $\mathcal{F}_T^V$ measurable random variable,
$$\tilde{E}_x[\la \ora{Y}_T, \phi \ra \varphi  ] = \int_{\mathbb{R}} p_0(x) \tilde{E}_x[ \phi(U_T) Q_T \varphi ] dx.$$
It then follows from the fact that $\tilde{P}_x(\cdot | \mathcal{F}_T^V) = \tilde{P}(\cdot | \mathcal{F}_T^V)$, and definition of $\tilde{P}$
$$\tilde{E}[\la \ora{Y}_T, \phi \ra \varphi] = \tilde{E}[\phi(U_T) Q_T \varphi], $$
as required in the theorem.
$\Box$

  {\bf Remark.}  From \eqref{eq:pro-transform}, we can see that 
$$
E\left[\phi(U_T)\big| \mathcal{F}_T^V\right]=\frac{\la\ora{Y}_T, Y_T \ra}{\tilde{\E}\left[ Q_t \big| \mathcal{F}_t^V \right]}
$$
Thus the solution $\ora{Y}_T$ of the FBDSDE \eqref{Forward:FBDSDEs} indeed provides an unnormalized solution for the optimal filter problem.  

\section{Closing Remarks}\label{sec:close} In this paper, we derived a Feymann-Kac type BDSDE formula for optimal filter problems and its adjoint. Then we show that the adjoint provides a unnormalized  solution for the optimal filter problem (BSDE filter).  As our preliminary work has shown, the BSDE filter has the potential to solve the optimal filter problem with more accuracy and less complexity than  traditional filter methods.

%

\end{document}